\documentclass[11pt,draftcls]{IEEEtran}{\onecolumn}
\usepackage{amsthm}
\usepackage{amsmath}
\usepackage{slashbox}
\usepackage{graphicx}
\usepackage{epstopdf}
\usepackage{epsfig,graphics,subfigure,psfrag,amsmath,amssymb}
\usepackage[nooneline,flushleft]{caption2}
\usepackage{cite}
\usepackage{float}
\usepackage{tocvsec2}
\usepackage{color}
\usepackage{verbatim}
\usepackage{algorithm, algorithmic}
\usepackage{mathtools}
\usepackage{hhline}
\usepackage{enumerate}
\usepackage{arydshln}

\newtheorem{thm}{Theorem}
\newtheorem{lem}{Lemma}

\newtheorem{rem}{Remark}

\newtheorem{assump}{Assumption}

\makeatletter

\newcommand{\Rmnum}[1]{\expandafter\@slowromancap\romannumeral #1@}
\makeatother

\DeclareMathOperator{\sgn}{sgn}
\DeclareMathOperator{\diag}{diag}

\begin{document}
\title{Decentralized Sparse Multitask RLS over Networks}
\author{Xuanyu Cao and K. J. Ray Liu, \emph{Fellow, IEEE}\\
Email: \{apogne, kjrliu\}@umd.edu\\
Department of Electrical and Computer Engineering, University of Maryland, College Park, MD\\}
\maketitle
\begin{abstract}
Distributed adaptive signal processing has attracted much attention in the recent decade owing to its effectiveness in many decentralized real-time applications in networked systems. Because many natural signals are highly sparse with most entries equal to zero, several decentralized sparse adaptive algorithms have been proposed recently. Most of them is focused on the single task estimation problems, in which all nodes receive data associated with the same unknown vector and collaborate to estimate it. However, many applications are inherently multitask oriented and each node has its own unknown vector different from others'. The related multitask estimation problem benefits from collaborations among the nodes as neighbor nodes usually share analogous properties and thus similar unknown vectors. In this work, we study the distributed sparse multitask recursive least squares (RLS) problem over networks. We first propose a decentralized online alternating direction method of multipliers (ADMM) algorithm for the formulated RLS problem. The algorithm is simplified for easy implementation with closed-form computations in each iteration and low storage requirements. Moreover, to further reduce the complexity, we present a decentralized online subgradient method with low computational overhead. We theoretically analyze the convergence behavior of the proposed subgradient method and derive an error bound related to the network topology and algorithm parameters. The effectiveness of the proposed algorithms is corroborated by numerical simulations and an accuracy-complexity tradeoff between the proposed two algorithms is highlighted.
\end{abstract}


\section{Introduction}

In the last decade, distributed adaptive signal processing has emerged as a vital topic because of the vast applications in need of decentralized real-time data processing over networked systems. For multi-agent networks, distributed adaptive algorithms only rely on local information exchange, i.e., information exchange among neighbor nodes, to estimate the unknowns. This trait endows distributed adaptive algorithms with low communication overhead, robustness to node/link failures and scalability to large networks. In the literature, the centralized least mean squares (LMS) and recursive least squares (RLS) \cite{Haykin:1996:AFT:230061} have been extended to their decentralized counterparts \cite{sayed2014adaptive,jiang2013distributed} to deal with estimation problems over networks. Furthermore, many natural signals are inherently sparse with most entries equal to zero such as the image signals and audio signals in \cite{duarte2008single,griffin2011single,donoho2006compressed,candes2006robust}. Sparsity of signals are particularly conspicuous in the era of big data: for many applications, redundant input features (e.g., a person's salary, education, height, gender, etc.) are collected to be fed into a learning system to predict a desired output (e.g., whether a person will resign his/her job). Most input features are unrelated to the output so that the weight vector between the input vector and the output is highly sparse. As such, several sparse adaptive algorithms have been proposed such as the sparse LMS in \cite{su2012performance,jin2010stochastic}, the sparse RLS in \cite{babadi2010sparls} and the distributed sparse RLS in \cite{liu2014distributed}.

Most of the decentralized sparse adaptive algorithms are focused on the single task estimation problem, in which all nodes receive data associated with the same unknown vector and collaborate to estimate it. On the contrary, many applications are inherently multitask-oriented, i.e., each node has its own unknown vector different from others'. For instance, in a sensor network, each node may want to estimate an unknown vector related to its specific location and thus different nodes have different unknown vectors to be estimated. In fact, several decentralized multitask adaptive algorithms have been proposed in the literature including the multitask diffusion LMS in \cite{chen2014multitask}, its asynchronous version in \cite{nassif2014multitask} and its application in the study of tremor in Parkinson's disease \cite{monajemi2016informed}. In particular, a sparse multitask LMS algorithm is proposed in \cite{nassif2015proximal} to promote sparsity of the estimated multitask vectors.

To the best of our knowledge, all the existing distributed adaptive algorithms for multitask estimation problems are based on various decentralized versions of LMS. The RLS based sparse multitask estimation problems have not aroused much attention. It is well known that the RLS possesses much faster convergence speed than the LMS. Hence, the RLS is more suitable for applications in need of fast and accurate tracking of the unknowns than the LMS, especially when the devices are capable of dealing with computations of moderately high complexity (which is the case as the computational capability of devices is increasing drastically). This motivates us to study the decentralized sparse multitask RLS problem over networks. The main contributions of this paper are summarized as follows.

\begin{itemize}
\item A global networked RLS minimization problem is formulated. In accordance with the multitask nature of the estimation problem, each node has its own weight vector. Since neighbor nodes often share analogous properties and thus similar weight vectors, we add regularization term to penalize deviations of neighbors' weight vectors. To enforce sparsity of the weight vectors, we further introduce $l_1$ regularization.
\item A decentralized online alternating direction method of multipliers (ADMM) algorithm is proposed for the formulated sparse multitask RLS problem. The proposed ADMM algorithm is simplified so that each iteration consists of simple closed-form computations and each node only needs to store and update one $M\times M$ matrix and six $M$ dimensional vectors, where $M$ is the dimension of the weight vectors. We show that the gaps between the outputs of the proposed ADMM algorithm and the optimal points of the formulated RLS problems converge to zero.
\item To overcome the relatively high computational cost of the proposed ADMM algorithm, we further present a decentralized online subgradient method, which enjoys lower computational complexity. We theoretically analyze its convergence behaviors and show that the tracking error of the weight vectors is upper bounded by some constant related to the network topology and algorithm parameters.
\item Numerical simulations are conducted to corroborate the effectiveness of the proposed algorithms. Their advantages over the single task sparse RLS algorithm in \cite{liu2014distributed} are highlighted. We also observe an accuracy-complexity tradeoff between the proposed two algorithms.
\end{itemize}

The roadmap of the remaining part of this paper is as follows. In Section \Rmnum{2}, the sparse multitask RLS problem is formally formulated. In Section \Rmnum{3}, we propose and simplify a decentralized online ADMM algorithm for the formulated RLS problem. In Section \Rmnum{4}, we propose a decentralized online subgradient method for the formulated problem in order to reduce computational complexity. In Section \Rmnum{5}, numerical simulations are conducted. In Section \Rmnum{6}, we conclude this work.

\section{The Statement of the Problem}

We consider a network of $N$ nodes and some edges between these nodes. We assume that the network is a simple graph, i.e., the network is undirected with no self-loop and there is at most one edge between any pair of nodes. Denote the set of neighbors of node $n$ (those who are linked with node $n$ by an edge) as $\Omega_n$. The network can be either connected or disconnected (there does not necessarily exist a path connecting every pair of nodes). Time is divided into discrete slots denoted as $t=1,2,...$. Each node $n$ has an unknown (slowly) time-variant $M$ dimensional weight vector $\widetilde{\mathbf{w}}_n(t)\in\mathbb{R}^M$ to be estimated. The formulated network is therefore a multitask learning network since different nodes have different weight vectors, as opposed to the traditional single task learning network \cite{sayed2014adaptive}, which is usually transformed into a consensus optimization problem framework \cite{shi2014linear,ling2015dlm,shi2015extra}. Each node $n$ has access to a sequence of private measurements $\{d_n(t),\mathbf{u}_n(t)\}_{t=1,2,...}$, where $\mathbf{u}_n(t)\in\mathbb{R}^M$ is the input regressor at time $t$ and $d_n(t)\in\mathbb{R}$ is the output observation at time $t$. The measurement data are private in the sense that node $n$ has access only to its own measurement sequence but not others'. The data at node $n$ are assumed to conform to a linear regression model with (slowly) time-variant weight vector $\widetilde{\mathbf{w}}_n(t)$:
\begin{equation}
d_n(t)=\mathbf{u}_n^\mathsf{T}(t)\widetilde{\mathbf{w}}_n(t)+e_n(t),
\end{equation}
where $e_n(t)$ is the output measurement noise at time $t$. In multitask learning networks, the benefit of cooperation between nodes comes from the fact that neighboring nodes have \emph{similar} weight vectors \cite{chen2014multitask}, where similarity is embodied by some specific distance measures. By incorporating terms promoting similarity between neighbors and enforcing cooperation in the network, an estimator may achieve potentially higher performance than its non-cooperative counterpart.

Moreover, many signals in practice are highly sparse, i.e., most entries in the signal are equal to zero, with examples encompassing image signals, audio signals, etc. \cite{duarte2008single,griffin2011single,donoho2006compressed,candes2006robust}. The sparsity of signals is especially conspicuous in today's big data era because redundant data are collected as input features among which most are unrelated to the targeted output, leading to sparsity of the corresponding weight vectors. Furthermore, as per convention in adaptive algorithms \cite{Haykin:1996:AFT:230061}, we assume that the weight vectors $\widetilde{\mathbf{w}}_n(t)$ varies with time very slowly. This suggests that past data are of great merit to estimate the current weight vector, which justifies the advantage of the RLS (studied in this paper) over the LMS (studied in all existing works on multitask estimation \cite{chen2014multitask,nassif2014multitask,nassif2015proximal,monajemi2016informed}) in terms of convergence speed.

In all, we propose an RLS based estimator to track the unknown weight vectors $\{\widetilde{\mathbf{w}}_n(t)\}_{n=1,2,...,N}$ while enforcing similarity between neighbors' weight vectors and sparsity of all weight vectors. The estimator at time $T$ is the optimal solution of the following optimization problem:
\begin{eqnarray}\label{RLS}
\text{Minimize}_{\mathbf{w}_1,...,\mathbf{w}_N}~~\sum_{n=1}^N\sum_{t=1}^T\lambda^{T-t}\left(d_n(t)-\mathbf{u}_n^\mathsf{T}(t)\mathbf{w}_n\right)^2+\beta\sum_{n=1}^N\sum_{m\in\Omega_n}\|\mathbf{w}_n-\mathbf{w}_m\|_2^2+\gamma\sum_{n=1}^N\|\mathbf{w}_n\|_1,
\end{eqnarray}
where $0<\lambda<1,\beta>0,\gamma>0$ are the forgetting factor of the RLS algorithm, regularization coefficient for similarity between neighbors' weight vectors and regularization coefficient for sparsity, respectively. If $\beta=\infty$, then problem \eqref{RLS} enforces consensus of weight vectors across nodes, and thus degenerates to the sparse RLS problem in \cite{liu2014distributed}. Note that the measurement data $\{\mathbf{u}_n(t),d_n(t)\}$ arrives in a sequential manner, which necessitates an online (real time) algorithm to solve \eqref{RLS} due to the prohibitive computation and storage cost of offline methods. Further note that the private measurement data are distributed among network nodes. Thus, a distributed algorithm for \eqref{RLS} is imperative as centralized algorithms are vulnerable to link failures and can incur large communication costs, not to mention the privacy concerns of the private data. Therefore, we are aimed at finding \emph{distributed online} algorithm for solving \eqref{RLS}. In the following two sections, we propose two different distributed online algorithms with complementary merits in accuracy and computational complexity.

\section{The Decentralized Online ADMM}

In this section, we propose an alternating direction method of multipliers (ADMM) based decentralized online algorithm for solving \eqref{RLS}. It is further simplified so that its iteration consists of simple closed-form computations and each node only needs to store and update one $M\times M$ matrix and six $M$ dimensional vectors. We show that the gaps between the outputs of the proposed ADMM algorithm and the optimal points of \eqref{RLS} converge to zero. Before the development of the algorithm, we first present some rudimentary knowledge of ADMM in the following subsection.

\subsection{Preliminaries of ADMM}
ADMM is an optimization framework widely applied to various signal processing applications, including wireless communications \cite{shen2012distributed}, power systems \cite{zhang2016admm} and multi-agent coordination \cite{chang2014proximal}. It enjoys fast convergence speed under mild technical conditions \cite{deng2016global} and is especially suitable for the development of distributed algorithms \cite{boyd2011distributed,bertsekas1989parallel}. ADMM solves problems of the following form:
\begin{eqnarray}\label{admm_prime}
\text{Minimize}_{\mathbf{x},\mathbf{z}} f(\mathbf{x})+g(\mathbf{z})~~\text{s.t.}~~\mathbf{Ax+Bz=c},
\end{eqnarray}
where $\mathbf{A}\in\mathbb{R}^{p\times n},B\in\mathbb{R}^{p\times m},c\in\mathbb{R}^p$ are constants and $\mathbf{x}\in\mathbb{R}^n,\mathbf{z}\in\mathbb{R}^m$ are optimization variables. $f:\mathbb{R}^n\mapsto\mathbb{R}$ and $g:\mathbb{R}^m\mapsto\mathbb{R}$ are two convex functions. The augmented Lagrangian can be formed as:
\begin{equation}
\mathfrak{L}_\rho(\mathbf{x,z,y})=f(\mathbf{x})+g(\mathbf{z})+\mathbf{y}^\mathsf{T}(\mathbf{Ax+Bz-c})+\frac{\rho}{2}\|\mathbf{Ax+Bz-c}\|_2^2,
\end{equation}
where $\mathbf{y}\in\mathbb{R}^p$ is the Lagrange multiplier and $\rho>0$ is some constant. The ADMM then iterates over the following three steps for $k\geq0$ (the iteration index):
\begin{eqnarray}
&&\mathbf{x}^{k+1}=\arg\min_\mathbf{x} \mathfrak{L}_\rho\left(\mathbf{x},\mathbf{z}^k,\mathbf{y}^k\right),\label{x_prime}\\
&&\mathbf{z}^{k+1}=\arg\min_\mathbf{z} \mathfrak{L}_\rho\left(\mathbf{x}^{k+1},\mathbf{z},\mathbf{y}^k\right),\label{z_prime}\\
&&\mathbf{y}^{k+1}=\mathbf{y}^{k}+\rho\left(\mathbf{Ax}^{k+1}+\mathbf{Bz}^{k+1}-\mathbf{c}\right).\label{multiplier_prime}
\end{eqnarray}
The ADMM is guaranteed to converge to the optimal point of \eqref{admm_prime} as long as $f$ and $g$ are convex \cite{boyd2011distributed,bertsekas1989parallel}. It is recently shown that global linear convergence can be ensured provided additional assumptions on problem \eqref{admm_prime} holds \cite{deng2016global}.

\subsection{Development of Decentralized Online ADMM for \eqref{RLS}}
To apply ADMM to \eqref{RLS}, we first transform it to the form of \eqref{admm_prime}. We introduce auxiliary variables $\mathbf{x}_n\in\mathbb{R}^M,n=1,...,N,$ and $\mathbf{v}_{n,i}\in\mathbb{R}^M,n=1,...N,i=1,...,|\Omega_n|$, where $|\cdot|$ denotes the cardinality of a set. Denote the index of $i$-th neighbor of node $n$ as $g(n,i)$. Thus, problem \eqref{RLS} can be equivalently transformed into the following problem:
\begin{eqnarray}
\begin{split}\label{admm}
&\text{Minimize}~~\sum_{n=1}^N\sum_{t=1}^T\lambda^{T-t}\left(d_n(t)-\mathbf{u}_n(t)^\mathsf{T}\mathbf{x}_n\right)^2\\
&~~~~~~~~~~~~+\beta\sum_{n=1}^N\left[|\Omega_n|\|\mathbf{x}_n\|_2^2-2\left(\sum_{i=1}^{|\Omega_n|}\mathbf{v}_{n,i}\right)^\mathsf{T}\mathbf{x}_n+\sum_{i=1}^{|\Omega_n|}\|\mathbf{v}_{n,i}\|_2^2\right]+\gamma\sum_{n=1}^N\|\mathbf{w}_n\|_1\\
&\text{s.t.}~~\mathbf{x}_n=\mathbf{w}_n,n=1,...,N,\\
&~~~~~~\mathbf{v}_{n,i}=\mathbf{w}_{g(n,i)},n=1...,N,i=1...,|\Omega_n|,
\end{split}
\end{eqnarray}
where the optimization variables are $\mathbf{w}_n,\mathbf{x}_n,\mathbf{v}_{n,i},n=1,...,N,i=1,...,|\Omega_n|$. Note that optimization problem \eqref{admm} is in the form of \eqref{admm_prime} (regarding $\mathbf{x}_n$'s and $\mathbf{v}_{n,i}$'s as the variable $\mathbf{x}$ in \eqref{admm_prime} and $\mathbf{w}_n$'s as the variable $\mathbf{z}$ in \eqref{admm_prime}). Thus, we can apply ADMM to problem \eqref{admm}. Introducing Lagrange multiplier $\mathbf{y}_n\in\mathbb{R}^M,\mathbf{z}_{n,i}\in\mathbb{R}^M$, we can form the augmented Lagrangian of \eqref{admm} as follows:
\begin{eqnarray}
\begin{split}
&\mathfrak{L}_\rho(\{\mathbf{x}_n,\mathbf{v}_{n,i},\mathbf{w}_n,\mathbf{y}_n,\mathbf{z}_{n,i}\}_{n=1,...,N,i=1,...,|\Omega_n|})\\
&=\sum_{n=1}^N\sum_{t=1}^T\lambda^{T-t}\left(d_n(t)-\mathbf{u}_n(t)^\mathsf{T}\mathbf{x}_n\right)^2+\beta\sum_{n=1}^N\left[|\mathcal{N}_n|\|\mathbf{x}_n\|_2^2-2\left(\sum_{i=1}^{|\Omega_n|}\mathbf{v}_{n,i}\right)^\mathsf{T}\mathbf{x}_n+\sum_{i=1}^{|\Omega_n|}\|\mathbf{v}_{n,i}\|_2^2\right]\\
&~~~+\gamma\sum_{n=1}^N\|\mathbf{w}_n\|_1+\sum_{n=1}^N\mathbf{y}_n^\mathsf{T}(\mathbf{x}_n-\mathbf{w}_n)+\sum_{n=1}^N\sum_{i=1}^{|\Omega_n|}\mathbf{z}_{n,i}^\mathsf{T}(\mathbf{v}_{n,i}-\mathbf{w}_{g_{n,i}})\\
&~~~+\frac{\rho}{2}\sum_{n=1}^N\|\mathbf{x}_n-\mathbf{w}_n\|_2^2+\frac{\rho}{2}\sum_{n=1}^N\sum_{i=1}^{|\Omega_n|}\|\mathbf{v}_{n,i}-\mathbf{w}_{g(n,i)}\|_2^2
\end{split}
\end{eqnarray}
In the following, for ease of notation, we use $\mathbf{x}$ to represent all the $\{\mathbf{x}_n\}$ and similarly for $\mathbf{v,w,y,z}$. We apply the ADMM updates \eqref{x_prime}, \eqref{z_prime} and \eqref{multiplier_prime} to problem \eqref{admm} as follows:
\begin{eqnarray}
&&\left\{\mathbf{x}^{k+1},\mathbf{v}^{k+1}\right\}=\arg\min_{\mathbf{x,v}} \mathfrak{L}_\rho\left(\mathbf{x},\mathbf{v},\mathbf{w}^k,\mathbf{y}^k,\mathbf{z}^k\right),\label{xv}\\
&&\mathbf{w}^{k+1}=\arg\min_{\mathbf{w}} \mathfrak{L}_\rho\left(\mathbf{x}^{k+1},\mathbf{v}^{k+1},\mathbf{w},\mathbf{y}^k,\mathbf{z}^k\right),\label{w}\\
&&\mathbf{y}_n^{k+1}=\mathbf{y}_n^k+\rho\left(\mathbf{x}_n^{k+1}-\mathbf{w}_n^{k+1}\right),\label{y}\\
&&\mathbf{z}_{n,i}^{k+1}=\mathbf{z}_{n,i}^k+\rho\left(\mathbf{v}_{n,i}^{k+1}-\mathbf{w}_{g(n,i)}^{k+1}\right).\label{z}
\end{eqnarray}
In the following, we detail how to implement the updates of the primal variables, i.e., \eqref{xv} and \eqref{w}, in a distributed and online fashion.

\subsubsection{Updating $\mathbf{x}$ and $\mathbf{v}$}
The update of $\mathbf{x}$ and $\mathbf{v}$ in \eqref{xv} can be decomposed across nodes. For each node $n$, the subproblem is:
\begin{eqnarray}
\begin{split}\label{xv1}
&\left\{\mathbf{x}_n^{k+1},\left\{\mathbf{v}_{n,i}^{k+1}\right\}_{i=1,...,|\Omega_n|}\right\}\\
&=\arg\min_{\mathbf{x}_n,\left\{\mathbf{v}_{n,i}\right\}_{i=1,...,|\Omega_n|}}~\Bigg\{\sum_{t=1}^T\lambda^{T-t}\left(d_n(t)-\mathbf{u}_n(t)^\mathsf{T}\mathbf{x}_n\right)^2\\
&+\beta\left[|\Omega_n|\|\mathbf{x}_n\|_2^2-2\left(\sum_{i=1}^{|\Omega_n|}\mathbf{v}_{n,i}\right)^\mathsf{T}\mathbf{x}_n+\sum_{i=1}^{|\Omega_n|}\|\mathbf{v}_{n,i}\|_2^2\right]\\
&+\mathbf{y}_n^{k\mathsf{T}}\mathbf{x}_n+\sum_{i=1}^{|\Omega_n|}\mathbf{z}_{n,i}^{k\mathsf{T}}\mathbf{v}_{n,i}+\frac{\rho}{2}\left\|\mathbf{x}_n-\mathbf{w}_n^k\right\|_2^2+\frac{\rho}{2}\sum_{i=1}^{|\Omega_n|}\left\|\mathbf{v}_{n,i}-\mathbf{w}_{g(n,i)}^k\right\|_2^2\Bigg\}.
\end{split}
\end{eqnarray}
Define the data dependent input correlation matrix and input-output cross correlation vector of node $n$ at time $T$ to be:
\begin{eqnarray}
&&\mathbf{R}_n(T)=\sum_{t=1}^T\lambda^{T-t}\mathbf{u}_n(T)\mathbf{u}_n(T)^\mathsf{T},\label{correlation}\\
&&\mathbf{p}_n(T)=\sum_{t=1}^T\lambda^{T-t}d_n(t)\mathbf{u}_n(t).\label{cross-correlation}
\end{eqnarray}
Note that the objective function in \eqref{xv1} is a convex quadratic function. Hence, the necessary and sufficient condition for optimality of problem \eqref{xv1} is that the gradient of the objective function vanishes. The gradient of the objective function, which is denoted as $J_n^k(T)$, with respect to $\mathbf{x}_n$ and $\mathbf{v}_{n,i}$ can be computed as follows:
\begin{eqnarray}
&&\nabla_{\mathbf{x}_n}J_n^k(T)=(2\mathbf{R_n(T)}+2\beta|\Omega_n|\mathbf{I}+\rho\mathbf{I})\mathbf{x}_n-2\beta\sum_{i=1}^{|\Omega_n|}\mathbf{v}_{n,i}-2\mathbf{p}_n(T)+\mathbf{y}_n^k-\rho\mathbf{w}_n^k,\\
&&\nabla_{\mathbf{v}_{n,i}}J_n^k(T)=-2\beta\mathbf{x}_n+(2\beta+\rho)\mathbf{v}_{n,i}+\mathbf{z}_{n,i}^k-\rho\mathbf{w}_{g(n,i)}^k.
\end{eqnarray}
Letting the gradients with respect to $\mathbf{x}_n$ and $\mathbf{v}_{n,i}$ be zero, we rewrite the update in \eqref{xv1} as:
\begin{eqnarray}\label{xv_matrix}
\begin{split}
&\left[
\begin{array}{c;{2pt/2pt}cccc}
2\mathbf{R}_n(T)+2\beta|\Omega_n|\mathbf{I}+\rho\mathbf{I} & -2\beta\mathbf{I} & -2\beta\mathbf{I} & \cdots & -2\beta\mathbf{I}\\ \hdashline[2pt/2pt]
-2\beta\mathbf{I} & (2\beta+\rho)\mathbf{I} &  &  & \\
-2\beta\mathbf{I} &  & (2\beta+\rho)\mathbf{I} & & \text{\huge0}\\
\vdots & & & \ddots &\\
-2\beta\mathbf{I} & &\text{\huge0} & & (2\beta+\rho)\mathbf{I}
\end{array}
\right]
\left[
\begin{array}{c}
\mathbf{x}_n^{k+1}\\ \hdashline[2pt/2pt]
\mathbf{v}_{n,1}^{k+1}\\
\mathbf{v}_{n,2}^{k+1}\\
\vdots\\
\mathbf{v}_{n,|\Omega_n|}^{k+1}
\end{array}
\right]
\\&=\left[
\begin{array}{c}
2\mathbf{p}_n(T)-\mathbf{y}_n^k+\rho\mathbf{w}_n^k\\ \hdashline[2pt/2pt]
-\mathbf{z}_{n,1}^k+\rho\mathbf{w}_{g(n,1)}^k\\
-\mathbf{z}_{n,2}^k+\rho\mathbf{w}_{g(n,2)}^k\\
\vdots\\
-\mathbf{z}_{n,|\Omega_n|}^k+\rho\mathbf{w}_{g(n,|\Omega_n|)}^k\\
\end{array}
\right]
\end{split}
\end{eqnarray}
To inverse the matrix in \eqref{xv_matrix}, we need to use the following matrix inversion lemma.
\begin{lem}
For arbitrary matrices $\mathbf{A}\in\mathbb{R}^{m\times m},\mathbf{B}\in\mathbb{R}^{m\times n},\mathbf{C}\in\mathbb{R}^{n\times m},\mathbf{D}\in\mathbb{R}^{n\times n}$ such that all the matrix inversions at the R.H.S. of \eqref{inversion_lemma} exist, we have:
\begin{eqnarray}\label{inversion_lemma}
\left[
\begin{array}{cc}
\mathbf{A} & \mathbf{B}\\
\mathbf{C} & \mathbf{D}
\end{array}
\right]^{-1}
=\left[
\begin{array}{cc}
\left(\mathbf{A}-\mathbf{B}\mathbf{D}^{-1}\mathbf{C}\right)^{-1} & -\left(\mathbf{A}-\mathbf{B}\mathbf{D}^{-1}\mathbf{C}\right)^{-1}\mathbf{BD}^{-1}\\
-\mathbf{D}^{-1}\mathbf{C}\left(\mathbf{A}-\mathbf{B}\mathbf{D}^{-1}\mathbf{C}\right)^{-1}
 & \mathbf{D}^{-1}+\mathbf{D}^{-1}\mathbf{C}\left(\mathbf{A}-\mathbf{B}\mathbf{D}^{-1}\mathbf{C}\right)^{-1}\mathbf{BD}^{-1}.
\end{array}
\right]
\end{eqnarray}
\end{lem}
Define a new matrix:
\begin{equation}\label{F}
\mathbf{F}_n(T)=\left[2\mathbf{R}_n(T)+\left(\rho+\frac{2\beta\rho|\Omega_n|}{2\beta+\rho}\right)\mathbf{I}\right]^{-1}.
\end{equation}
By invoking the matrix inversion lemma \eqref{inversion_lemma}, we can solve for the update \eqref{xv_matrix} in closed form:
\begin{align}
&\mathbf{x}_n^{k+1}=\mathbf{F}_n(T)\left(2\mathbf{p}_n(T)-\mathbf{y}_n^k+\rho\mathbf{w}_n^k\right)+\frac{2\beta}{2\beta+\rho}\mathbf{F}_n(T)\sum_{i=1}^{|\Omega_n|}\left(-\mathbf{z}_{n,i}^k+\rho\mathbf{w}_{g(n,i)}^k\right),\label{xx}\\
&\mathbf{v}_{n,i}^{k+1}=\frac{2\beta}{2\beta+\rho}\mathbf{F}_n(T)\left(2\mathbf{p}_n(T)-\mathbf{y}_n^k+\rho\mathbf{w}_n^k\right)+\frac{1}{2\beta+\rho}\left(-\mathbf{z}_{n,i}^k+\rho\mathbf{w}_{g(n,i)}^k\right)\nonumber\\
&~~~~~~~~~+\left(\frac{2\beta}{2\beta+\rho}\right)^2\mathbf{F}_n(T)\sum_{j=1}^{|\Omega_n|}\left(-\mathbf{z}_{n,j}^k+\rho\mathbf{w}_{g(n,j)}^k\right).\label{vv}
\end{align}

\subsubsection{Updating $\mathbf{w}$}
We note that the update for $\mathbf{w}$ in \eqref{w} can be decomposed not only across nodes but also across each entry of the vector $\mathbf{w}_n$. For each node $n$, the $l$-th entry of $\mathbf{w}_n$ can be updated as follows:
\begin{align}
&w_n^{k+1}(l)=\arg\min_{w_n(l)}\Bigg\{\gamma|w_n(l)|-y_n^k(l)w_n(l)-\left(\sum_{g(m,i)=n}z_{m,i}^k(l)\right)w_n(l)+\frac{\rho}{2}\left[w_n(l)-x_n^{k+1}(l)\right]^2\nonumber\\
&~~~~~~~~~~~~+\frac{\rho}{2}\sum_{g(m,i)=n}\left[w_n(l)-v_{m,i}^{k+1}(l)\right]^2\Bigg\}\\
&=\arg\min_{w_n(l)}\Bigg\{\gamma|w_n(l)|+\frac{\rho}{2}(1+|\Omega_n|)\Bigg[w_n(l)-\frac{1}{\rho(1+|\Omega_n|)}\Bigg(y_n^k(l)+\rho x_n^{k+1}(l)\nonumber\\
&~~~~~~~~~~~~+\sum_{g(m,i)=n}\left(z_{m,i}^k(l)+v_{m,i}^{k+1}(l)\right)\Bigg)\Bigg]^2\Bigg\}\\
&=\mathcal{S}_\frac{\gamma}{\rho(1+|\Omega_n|)}\left(\frac{1}{\rho(1+|\Omega_n|)}\Bigg(y_n^k(l)+\rho x_n^{k+1}(l)+\sum_{g(m,i)=n}\left(z_{m,i}^k(l)+v_{m,i}^{k+1}(l)\right)\Bigg)\right),\label{soft}
\end{align}
where the soft-threshold function $\mathcal{S}$ is defined for $a\in\mathbb{R},\kappa>0$ as follows:
\begin{align}
\mathcal{S}_\kappa(a)=
\begin{cases}
a-\kappa,\text{~~if~~}a>\kappa,\\
0,\text{~~if~~}|a|\leq\kappa,\\
a+\kappa,\text{~~if~~}a<\kappa.
\end{cases}
\end{align}
In \eqref{soft}, we have made use of the following fact.
\begin{lem}
For any $\lambda>0,\rho>0,v\in\mathbb{R}$, we have:
\begin{align}
\mathcal{S}_\frac{\lambda}{\rho}(v)=\arg\min_x~\left(\lambda |x|+\frac{\rho}{2}(x-v)^2\right).
\end{align}
\end{lem}
Once we extend the definition of $\mathcal{S}$ to vectors in a entrywise way, we can write the update for $\mathbf{w}_n$ compactly as:
\begin{align}
\mathbf{w}_n^{k+1}=\mathcal{S}_\frac{\gamma}{\rho(1+|\Omega_n|)}\left(\frac{1}{\rho(1+|\Omega_n|)}\Bigg(\mathbf{y}_n^k+\rho \mathbf{x}_n^{k+1}+\sum_{g(m,i)=n}\left(\mathbf{z}_{m,i}^k+\mathbf{v}_{m,i}^{k+1}\right)\Bigg)\right).\label{ww}
\end{align}

\subsubsection{Online Algorithm with Varying $T$}
So far, the derived ADMM algorithm is only suitable for one particular time shot $T$. Since it takes iterations for ADMM to converge to the optimal point, for each time $T$, we ought to run multiple rounds of ADMM iterations $k=1,...,K$ for some sufficiently large $K$. After the ADMM has converged for this particular time $T$, we update the data related quantities ($\mathbf{R}_n(T)$ and $\mathbf{p}_n(T)$) and move to the next time slot. However, since the underlying weight vectors are varying across time (i.e., the underlying linear system is non-stationary), it is meaningless to estimate the weight vectors very accurately for every time slot. Thus, in the following, we choose $K=1$, i.e., only one iteration of ADMM update is executed in each time slot. This is inspired by many existing adaptive algorithms such as the LMS algorithm, where only one step of gradient descent is performed at each time slot \cite{Haykin:1996:AFT:230061}. As such, we replace $k$ with $T-1$ in the previously derived updates \eqref{xx}, \eqref{vv}, \eqref{ww} and get updates that are suitable for varying time $T$:
\begin{align}
&\mathbf{x}_n(T)=\mathbf{F}_n(T)\left(2\mathbf{p}_n(T)-\mathbf{y}_n(T-1)+\rho\mathbf{w}_n(T-1)\right)\nonumber\\
&~~~~~~~~~+\frac{2\beta}{2\beta+\rho}\mathbf{F}_n(T)\sum_{i=1}^{|\Omega_n|}\left(-\mathbf{z}_{n,i}(T-1)+\rho\mathbf{w}_{g(n,i)}(T-1)\right),\label{xxx}\\
&\mathbf{v}_{n,i}(T)=\frac{2\beta}{2\beta+\rho}\mathbf{F}_n(T)\left(2\mathbf{p}_n(T)-\mathbf{y}_n(T-1)+\rho\mathbf{w}_n(T-1)\right)\nonumber\\
&~~~~~~~~~+\frac{1}{2\beta+\rho}\left(-\mathbf{z}_{n,i}(T-1)+\rho\mathbf{w}_{g(n,i)}(T-1)\right)\nonumber\\
&~~~~~~~~~+\left(\frac{2\beta}{2\beta+\rho}\right)^2\mathbf{F}_n(T)\sum_{j=1}^{|\Omega_n|}\left(-\mathbf{z}_{n,j}(T-1)+\rho\mathbf{w}_{g(n,j)}(T-1)\right),\label{vvv}\\
&\mathbf{w}_n(T)=\mathcal{S}_\frac{\gamma}{\rho(1+|\Omega_n|)}\left(\frac{1}{\rho(1+|\Omega_n|)}\Bigg(\mathbf{y}_n(T-1)+\rho \mathbf{x}_n(T)+\sum_{g(m,i)=n}\left(\mathbf{z}_{m,i}(T-1)+\mathbf{v}_{m,i}(T)\right)\Bigg)\right).\label{www}
\end{align}
Moreover, the updates \eqref{y} and \eqref{z} for dual variables can be rewritten as:
\begin{align}
\label{y_algo}&\mathbf{y}_n(T)=\mathbf{y}_n(T-1)+\rho\left(\mathbf{x}_n(T)-\mathbf{w}_n(T)\right),\\
&\mathbf{z}_{n,i}(T)=\mathbf{z}_{n,i}(T-1)+\rho\left(\mathbf{v}_{n,i}(T)-\mathbf{w}_{g(n,i)}(T)\right).\label{zzz}
\end{align}
The correlation matrices and cross-correlation vectors can be updated as follows:
\begin{align}
&\mathbf{R}_n(T+1)=\lambda\mathbf{R}_n(T)+\mathbf{u}_n(T+1)\mathbf{u}_n(T+1)^\mathsf{T},\label{R_update}\\
&\mathbf{p}_n(T+1)=\lambda\mathbf{p}_n(T)+d_n(T+1)\mathbf{u}_n(T+1).\label{p_update}
\end{align}
And $\mathbf{F}_n(T)$ is computed according to \eqref{F}.

\begin{rem}
The computation of $\mathbf{F}_n(T)$ in \eqref{F} necessitates inversion of an $M\times M$ matrix, which incurs a computational complexity of $\mathcal{O}(M^3)$ unless special structure is present. For the special case of $\lambda=1$ (which is suitable for time-invariant weight vectors), this burden can be alleviated as follows. According to \eqref{F}, \eqref{R_update} and the condition that $\lambda=1$, we have:
\begin{align}
\mathbf{F}_n(T)&=\left[2\left(\mathbf{R}_n(T-1)+\mathbf{u}_n(T)\mathbf{u}_n(T)^\mathsf{T}\right)+\left(\rho+\frac{2\beta\rho|\Omega_n|}{2\beta+\rho}\right)\mathbf{I}\right]^{-1}\\
&=\left[\mathbf{F}_n^{-1}(T-1)+2\mathbf{u}_n(T)\mathbf{u}_n(T)^\mathsf{T}\right]^{-1}\\
&=\mathbf{F}_n(T-1)-\frac{\mathbf{F}_n(T-1)\mathbf{u}_n(T)\mathbf{u}_n(T)^\mathsf{T}\mathbf{F}_n(T-1)}{\frac{1}{2}+\mathbf{u}_n(T)^\mathsf{T}\mathbf{F}_n(T-1)\mathbf{u}_n(T)}.\label{special}
\end{align}
However, in the general case where $\lambda<1$, the matrix inversion incurred by the computation of $\mathbf{F}_n(T)$ is inevitable, which is the most computationally intensive part of the proposed ADMM algorithm.
\end{rem}

\subsubsection{Simplification of the ADMM Updates}
So far, the ADMM updates involve primal variables $\{\mathbf{v}_{n,i}\}$ and dual variables $\{\mathbf{z}\}_{n,i}$. For each node $n$, $\{\mathbf{v}_{n,i}\}$ and $\{\mathbf{z}\}_{n,i}$ include $2|\Omega_n|$ $M$-dimensional vectors, which is costly to sustain in terms of communication and storage overhead, especially when the numbers of neighbors (degrees) are large. This motivates us to simplify the ADMM updates \eqref{xxx}-\eqref{p_update} so that the number of vectors at each node is independent of its degree. To this end, we first define the following auxiliary variables:
\begin{align}
&\mathbf{\underline{z}}_n(T)=\sum_{i=1}^{|\Omega_n|}\mathbf{z}_{n,i}(T),\\
&\mathbf{\overline{z}}_n(T)=\sum_{g(m,i)=n}\mathbf{z}_{m,i}(T),\\
&\mathbf{\underline{v}}_n(T)=\sum_{i=1}^{|\Omega_n|}\mathbf{v}_{n,i}(T),\\
&\mathbf{\overline{v}}_n(T)=\sum_{g(m,i)=n}\mathbf{v}_{m,i}(T),\\
\label{w_up}&\mathbf{\overline{w}}_n(T)=\sum_{m\in\Omega_n}\mathbf{w}_m(T),\\
\label{eta}&\boldsymbol{\eta}_n(T)=\mathbf{F}_n(T)(2\mathbf{p}_n(T)-\mathbf{y}_n(T-1)+\rho\mathbf{w}_n(T-1)),\\
\label{theta}&\boldsymbol{\theta}_n(T)=\mathbf{F}_n(T)(-\mathbf{\underline{z}}_n(T-1)+\rho\mathbf{\overline{w}}_n(T-1)),\\
\label{eta_up}&\boldsymbol{\overline{\eta}}_n(T)=\sum_{m\in\Omega_n}\boldsymbol{\eta}_m(T),\\
\label{theta_up}&\boldsymbol{\overline{\theta}}_n(T)=\sum_{m\in\Omega_n}\boldsymbol{\theta}_m(T).
\end{align}
Thus, the update for $\mathbf{x}$ in \eqref{xxx} can be rewritten as:
\begin{align}\label{x_algo}
\mathbf{x}_n(T)=\boldsymbol{\eta}_n(T)+\frac{2\beta}{2\beta+\rho}\boldsymbol{\theta}_n(T).
\end{align}
Using \eqref{vvv} yields the update for $\mathbf{\underline{v}}_n(T)$ and $\mathbf{\overline{v}}_n(T)$:
\begin{align}\label{v_under_algo}
\mathbf{\underline{v}}_n(T)&=\frac{2\beta|\Omega_n|}{2\beta+\rho}\boldsymbol{\eta}_n(T)+\left(\frac{2\beta}{2\beta+\rho}\right)^2|\Omega_n|\boldsymbol{\theta}_n(T)+\frac{1}{2\beta+\rho}(-\mathbf{\underline{z}}_n(T-1)+\rho\mathbf{\overline{w}}_n(T-1)).
\end{align}
\begin{align}\label{v_up_algo}
\mathbf{\overline{v}}_n(T)&=\frac{2\beta}{2\beta+\rho}\boldsymbol{\overline{\eta}}_n(T)+\left(\frac{2\beta}{2\beta+\rho}\right)^2\boldsymbol{\overline{\theta}}_n(T)+\frac{1}{2\beta+\rho}(-\mathbf{\overline{z}}_n(T-1)+\rho|\Omega_n|\mathbf{w}_n(T-1)).
\end{align}
The update for $\mathbf{w}_n(T)$ can be rewritten as:
\begin{align}\label{w_algo}
\mathbf{w}_n(T)=\mathcal{S}_\frac{\gamma}{\rho(1+|\Omega_n|)}\left(\frac{1}{\rho(1+|\Omega_n|)}(\mathbf{y}_n(T-1)+\rho\mathbf{x}_n(T)+\mathbf{\overline{z}}_n(T-1)+\rho\mathbf{\overline{v}}_n(T))\right).
\end{align}
Similarly, from \eqref{zzz}, we can spell out the updates for $\mathbf{\underline{z}}_n(T)$ and $\mathbf{\overline{z}}_n(T)$:
\begin{align}
\label{z_under_algo}&\mathbf{\underline{z}}_n(T)=\mathbf{\underline{z}}_n(T-1)+\rho(\mathbf{\underline{v}}_n(T)-\mathbf{\overline{w}}_n(T)),\\
\label{z_up_algo}&\mathbf{\overline{z}}_n(T)=\mathbf{\overline{z}}_n(T-1)+\rho(\mathbf{\overline{v}}_n(T)-|\Omega_n|\mathbf{w}_n(T)).
\end{align}
Now, we are ready to formally present the proposed decentralized online ADMM algorithm for solving \eqref{RLS}, which is summarized in Algorithm \ref{a1}. Notice that the algorithm is completely distributed: each node only needs to communicate with its neighbors. It is also online (real-time): each node only needs to store and update one $M\times M$ matrix $\mathbf{R}_n(T)$ and six $M$ dimensional vectors $\mathbf{p}_n(T),\mathbf{w}_n(T),\mathbf{\overline{w}}_n(T),\mathbf{y}_n(T),\mathbf{\underline{z}}_n(T),\mathbf{\overline{z}}_n(T)$. All other involved quantities in Algorithm \ref{a1} are intermediate and can be derived from these stored matrices and vectors.

\begin{algorithm}[!htbp]
\renewcommand{\algorithmicrequire}{\textbf{Inputs:} }
\renewcommand\algorithmicensure {\textbf{Outputs:} }
\caption{The proposed decentralized online ADMM algorithm}
\begin{algorithmic}[1]\label{a1}
\REQUIRE~~\\
Measurement data stream at each node $\{\mathbf{u}_n(t),d_n(t)\}$, $n=1,2,...,N$, $t=1,2,...$
\ENSURE~~\\
Estimates of the unknown weight vectors at each node $\{\mathbf{w}_n(T)\}$, $T=1,2,...$
\STATE Initialize $\mathbf{R}_n(0)=\mathbf{0}_{M\times M}$ and $\mathbf{p}_n(0)=\mathbf{w}_n(0)=\mathbf{\overline{w}}_n(0)=\mathbf{y}_n(0)=\mathbf{\underline{z}}_n(0)=\mathbf{\overline{z}}_n(0)=\mathbf{0}_M$. $T=0$.
\STATE \textbf{Repeat:}
\STATE $T\leftarrow T+1$.
\STATE Each node $n$ updates its correlation matrix and cross-correlation vector once receiving the new data $\mathbf{u}_n(T),d_n(T)$:
\begin{align}
&\mathbf{R}_n(T)=\lambda\mathbf{R}_n(T-1)+\mathbf{u}_n(T)\mathbf{u}_n(T)^\mathsf{T},\\
&\mathbf{p}_n(T)=\lambda\mathbf{p}_n(T-1)+d_n(t)\mathbf{u}_n(T).
\end{align}
\STATE Each node $n$ computes $\mathbf{F}_n(T)$ according to \eqref{F}.
\STATE Each node $n$ computes $\boldsymbol{\eta}_n(T)$ and $\boldsymbol{\theta}_n(T)$ according to \eqref{eta} and \eqref{theta} and then broadcasts its the results to its neighbors.
\STATE Each node $n$ receives $\boldsymbol{\eta}_m(T)$ and $\boldsymbol{\theta}_m(T)$ from its neighbors $m\in\Omega_n$ and forms $\boldsymbol{\overline{\eta}}_n(T)$ and $\boldsymbol{\overline{\theta}}_n(T)$ based on \eqref{eta_up} and \eqref{theta_up}.
\STATE Each node $n$ computes $\mathbf{x}_n(T),\mathbf{\underline{v}}_n(T),\mathbf{\overline{v}}_n(T)$ according to \eqref{x_algo}, \eqref{v_under_algo} and \eqref{v_up_algo}, respectively.
\STATE Each node $n$ computes $\mathbf{w}_n(T)$ based on \eqref{w_algo} and broadcasts the result to its neighbors.
\STATE Each node $n$ receives $\mathbf{w}_m(T)$ from its neighbors $m\in\Omega_n$ and form $\mathbf{\overline{w}}_n(T)$ according to \eqref{w_up}.
\STATE Each node $n$ updates $\mathbf{y}_n(T),\mathbf{\underline{z}}_n(T),\mathbf{\overline{z}}_n(T)$ according to \eqref{y_algo}, \eqref{z_under_algo} and \eqref{z_up_algo}.
\end{algorithmic}
\end{algorithm}

\subsection{Convergence of the Algorithm \ref{a1}}
In this subsection, we briefly discuss about the convergence of Algorithm \ref{a1} and show that the gap between its output, $\mathbf{w}_n(T)$, and the optimal point of problem \eqref{RLS}, which we denote as $\mathbf{w}_n^*(T)$, converges to zero. We make the following assumptions.
\begin{assump}
The true weight vector $\mathbf{\widetilde{w}}_n$ is time-invariant, i.e., the linear regression data model is $d_n(t)=\mathbf{u}_n(t)^\mathsf{T}\mathbf{\widetilde{w}}_n+e_n(t)$.
\end{assump}
\begin{assump}
For each node $n$, the input process $\{\mathbf{u}_n(t)\}_{t=1,2,...}$ is independent across time with time-invariant correlation matrix $\mathbf{R}_n=\mathbb{E}\left[\mathbf{u}_n(t)\mathbf{u}_n(t)^\mathsf{T}\right]$.
\end{assump}
\begin{assump}
For each node $n$, the noise process $\{e_n(t)\}_{t=1,2,...}$ has zero mean, i.e., $\mathbb{E}[e_n(t)]=0$ and is independent across time and independent from the input process $\{\mathbf{u}_n(t)\}$.
\end{assump}

Note that all of these assumptions are standard when analyzing the performance of adaptive algorithms in the literature \cite{Haykin:1996:AFT:230061}. From the definition of $\mathbf{R}_n(T)$ and $\mathbf{p}_n(T)$, we know that they are weighted sum of i.i.d. terms. According to the strong law of large numbers for weighted sums \cite{chow1973limiting,babadi2010sparls}, as $n\rightarrow\infty$, $\mathbf{R}_n(T)$ converges to $\lim_{T\rightarrow\infty}\mathbb{E}[\mathbf{R}_n(T)]=\frac{\mathbf{R}_n}{1-\lambda}$. Similarly, $\mathbf{p}_n(T)$ converges to $\lim_{T\rightarrow\infty}\mathbb{E}[\mathbf{p}_n(T)]=\frac{\mathbf{R}_n\mathbf{\widetilde{w}}_n}{1-\lambda}$. When $T\rightarrow\infty$, the optimization problem at time $T$, i.e., problem \eqref{RLS}, is to minimize (w.r.t. $\{\mathbf{w}_n\}$):
\begin{align}
&\sum_{n=1}^N\sum_{t=1}^T\lambda^{T-t}\left(\mathbf{w}_n^\mathsf{T}\mathbf{u}_n(t)\mathbf{u}_n(t)^\mathsf{T}\mathbf{w}_n-2d_n(t)\mathbf{u}_n(t)^\mathsf{T}\mathbf{w}_n\right)+\beta\sum_{n=1}^N\sum_{m\in\Omega_n}\|\mathbf{w}_n-\mathbf{w}_m\|_2^2+\gamma\sum_{n=1}^N\|\mathbf{w}_n\|_1\nonumber\\
&\approx\frac{1}{1-\lambda}\sum_{n=1}^N\left(\mathbf{w}_n^\mathsf{T}\mathbf{R}_n\mathbf{w}_n-2\mathbf{\widetilde{w}}_n^\textsf{T}\mathbf{R}_n\mathbf{w}_n\right)+\beta\sum_{n=1}^N\sum_{m\in\Omega_n}\|\mathbf{w}_n-\mathbf{w}_m\|_2^2+\gamma\sum_{n=1}^N\|\mathbf{w}_n\|_1.\label{convergence_admm}
\end{align}
Note that the R.H.S. of \eqref{convergence_admm} does not depend on $T$. Note that ADMM is guaranteed to converge to the optimal point for static convex optimization problem of the form \eqref{admm_prime} and the R.H.S. of \eqref{convergence_admm} can be transformed into the form of \eqref{admm_prime} as we do in Subsection III-B. So, the output of Algorithm \ref{a1}, $\mathbf{w}_n(t)$, converges to the minimum point of the R.H.S. of \eqref{convergence_admm}. Due to \eqref{convergence_admm} and the definition of $\mathbf{w}_n^*(T)$, we know that $\mathbf{w}_n^*(T)$ also converges to the minimum point of the R.H.S. of \eqref{convergence_admm}. Hence, the difference between the output of Algorithm \ref{a1}, i.e., $\mathbf{w}_n(t)$, and  the optimal point of \eqref{RLS}, i.e., $\mathbf{w}_n^*(T)$, converges to zero.

\section{The Decentralized Online Subgradient Method}

The implementation of the proposed Algorithm \ref{a1} necessitates an inversion of an $M\times M$ matrix at each time and each node, which may not be suitable for nodes with low computational capability. In fact, a relatively high computational overhead is a general drawback of dual domain methods (e.g., ADMM) in optimization theory \cite{ling2015dlm}. On the contrary, primal domain methods such as gradient descent method, though having relatively slow convergence speed, enjoys low computational complexity \cite{nedic2009distributed}. As such, in this section, we present a distributed online subgradient method for problem \eqref{RLS} to trade off convergence speed and accuracy for low computational complexity.

\subsection{Development of the Decentralized Online Subgradient Method}

Recall the optimization problem at time $T$, i.e., problem \eqref{RLS}. Denote the objective function of \eqref{RLS} as $H_T(\mathbf{w})$. We derive the subdifferential (the set of subgradients \cite{boyd2006subgradient}) of $H_T$ at $\mathbf{w}$ to be:
\begin{align}\label{subdifferential}
\partial H_T(\mathbf{w})=
\left[
\begin{array}{c}
2\mathbf{R}_1(T)\mathbf{w}_1-2\mathbf{p}_1(T)+2\beta\left(2|\Omega_1|\mathbf{w}_1-2\sum_{m\in\Omega_1}\mathbf{w}_m\right)+\gamma\sgn(\mathbf{w}_1)\\
\vdots\\
2\mathbf{R}_N(T)\mathbf{w}_N-2\mathbf{p}_N(T)+2\beta\left(2|\Omega_N|\mathbf{w}_N-2\sum_{m\in\Omega_N}\mathbf{w}_m\right)+\gamma\sgn(\mathbf{w}_N)
\end{array}
\right],
\end{align}
where the sign (set) function is defined as:
\begin{align}\label{definition_sgn}
\sgn(x)=
\begin{cases}
1,~~\text{if}~~x>0,\\
-1,~~\text{if}~~x<0,\\
[-1,1],~~\text{if}~~x=0.
\end{cases}
\end{align}
The extension of the $\sgn$ function to vectors is entrywise. The subgradient method is to simply use the iteration $\mathbf{w}(T)=\mathbf{w}(T-1)-\alpha \mathbf{g}$, where $\mathbf{g}\in H_T(\mathbf{w}(T-1))$ is any subgradient of $H_T$ at $\mathbf{w}(T-1)$ and $\alpha>0$ is the step size \cite{boyd2006subgradient}. This naturally leads to the following decentralized online update:
\begin{align}\label{subgradient_w}
&\mathbf{w}_n(T)=\mathbf{w}_n(T-1)-\alpha\Bigg[2\mathbf{R}_n(T)\mathbf{w}_n(T-1)-2\mathbf{p}_n(T)+4\beta\sum_{m\in\Omega_n}(\mathbf{w}_n(T-1)-\mathbf{w}_m(T-1))\nonumber\\
&~~~~~~~~~~~~~+\gamma\sgn(\mathbf{w}_n(T-1))\Bigg],
\end{align}
where $\sgn(0)$ is any number within the interval $[-1,1]$\footnote{There is a standard abuse of notation for the $\sgn$ function: in \eqref{subdifferential} and \eqref{definition_sgn}, $\sgn(0)$ is defined to be the interval $[-1,1]$ while in \eqref{subgradient_w}, $\sgn(0)$ is defined to be any arbitrary number within $[-1,1]$. In the following, the latter definition will be used.}. By introducing an auxiliary variable $\mathbf{\overline{w}}_n(T)$, we propose the decentralized online subgradient method for \eqref{RLS} in Algorithm \ref{a2}. We observe that Algorithm \ref{a2} is completely decentralized as every node only communicates with its neighbors. It is also online since each node only needs to store and update one $M\times M$ matrix and three $M$ dimensional vectors. More importantly, Algorithm \ref{a2} is free of any matrix inversion, which is a major burden of Algorithm \ref{a1}.

\begin{algorithm}[!htbp]
\renewcommand{\algorithmicrequire}{\textbf{Inputs:} }
\renewcommand\algorithmicensure {\textbf{Outputs:} }
\caption{The proposed decentralized online subgradient algorithm}
\begin{algorithmic}[1]\label{a2}
\REQUIRE~~\\
Measurement data stream at each node $\{\mathbf{u}_n(t),d_n(t)\}$, $n=1,2,...,N$, $t=1,2,...$
\ENSURE~~\\
Estimates of the unknown weight vectors at each node $\{\mathbf{w}_n(T)\}$, $T=1,2,...$
\STATE Initialize $\mathbf{R}_n(0)=\mathbf{0}_{M\times M}$, $\mathbf{p}_n(0)=\mathbf{w}_n(0)=\mathbf{\overline{w}}_n(0)=\mathbf{0}_M$, $T=0$.
\STATE \textbf{Repeat:}
\STATE $T\leftarrow T+1$.
\STATE Each node $n$ updates its correlation matrix and cross-correlation vector once receiving the new data $\mathbf{u}_n(T),d_n(T)$:
\begin{align}
&\mathbf{R}_n(T)=\lambda\mathbf{R}_n(T-1)+\mathbf{u}_n(T)\mathbf{u}_n(T)^\mathsf{T},\\
&\mathbf{p}_n(T)=\lambda\mathbf{p}_n(T-1)+d_n(t)\mathbf{u}_n(T).
\end{align}
\STATE Each node $n$ updates $\mathbf{w}_n(T)$:
\begin{align}
&\mathbf{w}_n(T)=\mathbf{w}_n(T-1)-\alpha\Bigg[2\mathbf{R}_n(T)\mathbf{w}_n(T-1)-2\mathbf{p}_n(T)+4\beta(|\Omega_n|\mathbf{w}_n(T-1)-\mathbf{\overline{w}}_n(T-1))\nonumber\\
&~~~~~~~~~~~~~+\gamma\sgn(\mathbf{w}_n(T-1))\Bigg].\label{update_subgradient}
\end{align}
\STATE Each node $n$ broadcasts its $\mathbf{w}_n(T)$ to its neighbors.
\STATE Each node $n$ receives $\mathbf{w}_m(T)$ from its neighbors $m\in\Omega_n$ and form $\mathbf{\overline{w}}_n(T)$ as follows:
\begin{align}
\mathbf{\overline{w}}_n(T)=\sum_{m\in\Omega_n}\mathbf{w}_m(T).
\end{align}
\end{algorithmic}
\end{algorithm}

\subsection{Convergence Analysis of Algorithm \ref{a2}}
In this subsection, we analyze the convergence behavior of Algorithm \ref{a2}. We still make the same assumptions as in the analysis of the proposed ADMM algorithm in Subsection III-C, i.e., Assumptions 1-3. As explained in Subsection III-C, for large $T$, we can approximate $\mathbf{R}_n(T)$ by $\frac{\mathbf{R}_n}{1-\lambda}$ and $\mathbf{p}_n(T)$ by $\frac{\mathbf{R}_n\mathbf{\widetilde{w}}_n}{1-\lambda}$. Define the error vector of node $n$ at time $T$ to be:
\begin{align}
\mathbf{f}_n(T)=\mathbf{w}_n(T)-\mathbf{\widetilde{w}}_n.\label{error}
\end{align}
Thus, substituting \eqref{error} into the definition of $\mathbf{\overline{w}}_n(T-1)$ yields:
\begin{align}
\mathbf{\overline{w}}_n(T-1)=\sum_{m\in\Omega_n}\left(\mathbf{\widetilde{w}}_m+\mathbf{f}_m(T-1)\right).
\end{align}
Hence, using \eqref{update_subgradient}, \eqref{error} and the approximation of $\mathbf{R}_n(T),\mathbf{p}_n(T)$, we can derive a recursive equation for the error vector:
\begin{align}
\mathbf{f}_n(T)&=\mathbf{w}_n(T-1)-\mathbf{\widetilde{w}}_n-\alpha\Bigg[2\mathbf{R}_n(T)\left(\mathbf{\widetilde{w}}_n+\mathbf{f}_n(T-1)\right)-2\mathbf{p}_n(T)\nonumber\\
&~~~+4\beta\left(|\Omega_n|\left(\mathbf{\widetilde{w}}_n+\mathbf{f}_n(T-1)\right)-\sum_{m\in\Omega_n}\left(\mathbf{\widetilde{w}}_m+\mathbf{f}_m(T-1)\right)\right)+\gamma\sgn(\mathbf{w}_n(T-1))\Bigg]\\
&\approx\mathbf{f}_n(T-1)-\alpha\Bigg[\frac{2\mathbf{R}_n}{1-\lambda}(\mathbf{\widetilde{w}}_n+\mathbf{f}_n(T-1))-\frac{2\mathbf{R}_n\mathbf{\widetilde{w}}_n}{1-\lambda}\nonumber\\
&~~~+4\beta\left(|\Omega_n|\mathbf{\widetilde{w}}_n-\sum_{m\in\Omega_n}\mathbf{\widetilde{w}}_m+|\Omega_n|\mathbf{f}_n(T-1)-\sum_{m\in\Omega_n}\mathbf{f}_m(T-1)\right)+\gamma\sgn(\mathbf{w}_n(T-1))\Bigg]\\
&=\left(\mathbf{I}-\frac{2\alpha\mathbf{R}_n}{1-\lambda}\right)\mathbf{f}_n(T-1)-4\alpha\beta\left(|\Omega_n|\mathbf{f}_n(T-1)-\sum_{m\in\Omega_n}\mathbf{f}_m(T-1)\right)\nonumber\\
&~~~-4\alpha\beta\left(|\Omega_n|\mathbf{\widetilde{w}}_n-\sum_{m\in\Omega_n}\mathbf{\widetilde{w}}_m\right)-\alpha\gamma\sgn(\mathbf{w}_n(T-1)).
\end{align}
Taking expectations yields:
\begin{align}
\label{expect}\mathbb{E}[\mathbf{f}_n(T)]=\left[(1-4\alpha\beta |\Omega_n|)\mathbf{I}-\frac{2\alpha\mathbf{R}_n}{1-\lambda}\right]\mathbb{E}[\mathbf{f}_n(T-1)]+4\alpha\beta\sum_{m\in\Omega_n}\mathbb{E}[\mathbf{f}_m(T-1)]+\mathbf{r}_n(T),
\end{align}
where $\mathbf{r}_n(T)$ is defined as:
\begin{align}\label{r_def}
\mathbf{r}_n(T)=-4\alpha\beta\left(|\Omega_n|\mathbf{\widetilde{w}}_n-\sum_{m\in\Omega_n}\mathbf{\widetilde{w}}_m\right)-\alpha\gamma\mathbb{E}[\sgn(\mathbf{w}_n(T-1))].
\end{align}
Denote the adjacency matrix of the network as $\mathbf{A}\in\mathbb{R}^{M\times M}$, i.e., $A(i,j)=1$ if $i,j$ are connected with an edge; otherwise $A(i,j)=0$. Define $\mathbf{C}=\mathbf{A}\otimes\mathbf{I}_{M\times M}\in\mathbb{R}^{MN\times MN}$ to be the block adjacency matrix, where $\otimes$ means Kronecker product. Define the matrix $\mathbf{B}\in\mathbb{R}^{MN\times MN}$ as follows:
\begin{align}
\mathbf{B}=
\left[
\begin{array}{lllc}
(1-4\alpha\beta|\Omega_1|)\mathbf{I}-\frac{2\alpha\mathbf{R}_1}{1-\lambda} & & &\\
& (1-4\alpha\beta|\Omega_2|)\mathbf{I}-\frac{2\alpha\mathbf{R}_2}{1-\lambda} & &\text{\huge0}\\
& & \ddots&\\
&\text{\huge0} & & (1-4\alpha\beta|\Omega_N|)\mathbf{I}-\frac{2\alpha\mathbf{R}_N}{1-\lambda}
\end{array}
\right]
\end{align}
Furthermore, stack $\mathbf{f}_n(T)$ and $\mathbf{r}_n(T)$ into long vectors respectively:
\begin{align}
\mathbf{f}(T)=
\left[
\begin{array}{c}
\mathbf{f}_1(T)\\
\vdots\\
\mathbf{f}_N(T)
\end{array}
\right]\in\mathbb{R}^{NM},
\mathbf{r}(T)=
\left[
\begin{array}{c}
\mathbf{r}_1(T)\\
\vdots\\
\mathbf{r}_N(T)
\end{array}
\right]\in\mathbb{R}^{NM}.
\end{align}

Hence, \eqref{expect} can be written in a more compact form:
\begin{align}\label{expect_big}
\mathbb{E}[\mathbf{f}(T)]=(\mathbf{B}+4\alpha\beta\mathbf{C})\mathbb{E}[\mathbf{f}(T-1)]+\mathbf{r}(T).
\end{align}
Define $\boldsymbol{\Phi}=\frac{1}{\alpha}(\mathbf{I}-\mathbf{B})$. Since $\boldsymbol{\Phi}-4\beta\mathbf{C}$ is symmetric, we can perform eigendecomposition for it, i.e., $\boldsymbol{\Phi}-4\beta\mathbf{C}=\mathbf{Q\Lambda Q}^\mathsf{T}$ for some orthogonal matrix $\mathbf{Q}\in\mathbb{R}^{MN\times MN}$ and diagonal matrix $\boldsymbol{\Lambda}=\diag(\delta_1,...,\delta_{MN})$, where $\delta_i$'s are the eigenvalues of $\boldsymbol{\Phi}-4\beta\mathbf{C}$. Therefore, $\mathbf{B}+4\alpha\beta\mathbf{C}=\mathbf{I}-\alpha(\mathbf{\Phi}-4\beta\mathbf{C})=\mathbf{Q\Sigma Q}^\mathsf{T}$, where we define $\mathbf{\Sigma}=\diag(1-\alpha\delta_1,...,1-\alpha\delta_{MN})$. Hence, \eqref{expect_big} can be rewritten as:
\begin{align}
\mathbf{Q}^\mathsf{T}\mathbb{E}[\mathbf{f}(T)]=\mathbf{\Sigma Q^\mathsf{T}}\mathbb{E}[\mathbf{f}(T-1)]+\mathbf{Q}^\mathsf{T}\mathbf{r}(T),
\end{align}
or
\begin{align}\label{epsilon_recursive}
\boldsymbol{\epsilon}(T)=\boldsymbol{\Sigma\epsilon}(T-1)+\mathbf{\widetilde{r}}(T),
\end{align}
where we define $\boldsymbol{\epsilon}(T)=\mathbf{Q}^\mathsf{T}\mathbb{E}[\mathbf{f}(T)]$ and $\mathbf{\widetilde{r}}(T)=\mathbf{Q}^\mathsf{T}\mathbf{r}(T)$. We first want to bound $\mathbf{\widetilde{r}}(T)$. To this end, from \eqref{r_def} we derive:
\begin{align}
\|\mathbf{r}_n(T)\|_2&\leq\alpha\mathbb{E}\left[4\beta|\Omega_n|\|\mathbf{\widetilde{w}}_n\|_2+4\beta\sum_{m\in\Omega_n}\|\mathbf{\widetilde{w}}_m\|_2+\gamma\|\sgn(\mathbf{w}_n(T-1))\|_2\right]\\
&\leq\alpha\left(8\beta\max_{n=1,...,N}|\Omega_n|\max_{n=1,...,N}\|\mathbf{\widetilde{w}}_n\|_2+\gamma\sqrt{M}\right).
\end{align}
Therefore,
\begin{align}\label{r_tilde_bound}
\|\mathbf{\widetilde{r}}(T)\|_2^2=\|\mathbf{Q}^\mathsf{T}\mathbf{r}(T)\|_2^2=\|\mathbf{r}(T)\|_2^2=\sum_{n=1}^N\|\mathbf{r}_n(T)\|_2^2\leq N\alpha^2\left(8\beta\max_{n=1,...,N}|\Omega_n|\max_{n=1,...,N}\|\mathbf{\widetilde{w}}_n\|_2+\gamma\sqrt{M}\right)^2.
\end{align}
Moreover, recursive applications of \eqref{epsilon_recursive} yields:
\begin{align}
\boldsymbol{\epsilon}(T)=\mathbf{\Sigma}^T\boldsymbol{\epsilon}(0)+\sum_{t=0}^{T-1}\boldsymbol{\Sigma}^t\mathbf{\widetilde{r}}(T-t),
\end{align}
where the superscript $T$ of $\mathbf{\Sigma}$ means power instead of transposition. So,
\begin{align}\label{key}
\|\boldsymbol{\epsilon}(T)\|_2\leq\|\mathbf{\Sigma}\|_2^T\|\boldsymbol{\epsilon}(0)\|_2+\sum_{t=0}^{T-1}\|\mathbf{\Sigma}\|_2^t\|\mathbf{\widetilde{r}}(T-t)\|_2,
\end{align}
where $\|\mathbf{\Sigma}\|_2$ means the spectral norm of $\mathbf{\Sigma}$, i.e., the maximum singular value of $\mathbf{\Sigma}$. Since $\mathbf{\Sigma}$ is a diagonal matrix, it is easy to calculate its spectral norm as $\|\mathbf{\Sigma}\|_2=\max_{i=1,...,NM}|1-\alpha\delta_i|$. Now, we make two additional assumptions in the following.
\begin{assump}
$\frac{2}{1-\lambda}\min_n\lambda_{\min}(\mathbf{R}_n)-4\beta(\max_{n}|\Omega_n|-\min_{n}|\Omega_n|)>0$.
\end{assump}
\begin{assump}
$\alpha<\frac{2}{\max_{i=1,...,NM}\delta_i}$.
\end{assump}
We note the following Weyl's theorem \cite{horn2012matrix}.
\begin{lem}
(Weyl) Given two symmetric matrices $\mathbf{A,B}\in\mathbb{R}^{n\times n}$. Arrange the eigenvalues in increasing order, i.e., $\lambda_1(\mathbf{X})\leq...\leq\lambda_n(\mathbf{X})$ for $\mathbf{X}\in\{\mathbf{A,B,A+B}\}$. Then for any $k\in\{1,...,n\}$, we have:
\begin{align}
\lambda_k(\mathbf{A})+\lambda_1(\mathbf{B})\leq\lambda_k(\mathbf{A+B})\leq\lambda_k(\mathbf{A})+\lambda_n(\mathbf{B}).
\end{align}
\end{lem}
Making use of Weyl's theorem, we obtain:
\begin{align}\label{eigen}
\lambda_{\min}(\mathbf{\Phi}-4\beta\mathbf{C})\geq\lambda_{\min}(\Phi)+\lambda_{\min}(-4\beta\mathbf{C})=\lambda_{\min}(\Phi)-4\beta\lambda_{\max}(\mathbf{C}).
\end{align}
Recall that $\mathbf{\Phi}=\diag(4\beta|\Omega_1|\mathbf{I}+\frac{2\mathbf{R}_1}{1-\lambda},...,4\beta|\Omega_N|\mathbf{I}+\frac{2\mathbf{R}_N}{1-\lambda})$. We have $\lambda_{\min}(\Phi)=\min_n\left\{4\beta|\Omega_n|+\frac{2}{1-\lambda}\lambda_{\min}(\mathbf{R}_n)\right\}$. Further noting that $\mathbf{C}=\mathbf{A}\otimes\mathbf{I}$ and defining $\|\mathbf{A}\|_1=\max_{1\leq j\leq N}\sum_{i=1}^N|A_{ij}|$ to be the maximum-column-sum norm of matrix $\mathbf{A}$, we have
\begin{align}
\lambda_{\max}(\mathbf{C})=\lambda_{\max}(\mathbf{A})\leq|\lambda_{\max}(\mathbf{A})|\leq\rho(\mathbf{A})\leq\|\mathbf{A}\|_1=\max_n|\Omega_n|,
\end{align}
where $\rho(\mathbf{A})=\max_{i=1,...,N}|\lambda_i(\mathbf{A})|$ is the spectral radius of $\mathbf{A}$. So, continuing from \eqref{eigen} we have:
\begin{align}
\lambda_{\min}(\mathbf{\Phi}-4\beta\mathbf{C})&\geq\min_n\left\{4\beta|\Omega_n|+\frac{2}{1-\lambda}\lambda_{\min}(\mathbf{R}_n)\right\}\\
&\geq\frac{2}{1-\lambda}\min_n\lambda_{\min}(\mathbf{R}_n)-4\beta(\max_{n}|\Omega_n|-\min_{n}|\Omega_n|)\\
&>0,
\end{align}
where the last step follows from Assumption 4. Since $\delta_i$'s are eigenvalues of $\mathbf{\Phi}-4\beta\mathbf{C}$, we thus have $\delta_i>0,i=1,...,MN$. This together with Assumption 5 implies that $\max_{i=1,...,NM}|1-\alpha\delta_i|<1$, i.e., $\|\mathbf{\Sigma}\|_2<1$. Hence, as $T\rightarrow\infty$, the first time of the R.H.S. of \eqref{key} converges to zero. Furthermore, since \eqref{r_tilde_bound} holds for any $T$, the second term of the R.H.S. of \eqref{key} is always bounded above by $\frac{\sqrt{N}\alpha}{1-\|\mathbf{\Sigma}\|_2}\left(8\beta\max_{n=1,...,N}|\Omega_n|\max_{n=1,...,N}\|\mathbf{\widetilde{w}}_n\|_2+\gamma\sqrt{M}\right)$. So, in all, from \eqref{key} and the fact that $\|\mathbb{E}[\mathbf{f}(T)]\|_2=\|\mathbf{Q}^\mathsf{T}\mathbb{E}[\mathbf{f}(T)]\|_2=\|\boldsymbol{\epsilon}(T)\|_2$, we have:
\begin{align}
\limsup_{T\rightarrow\infty}\|\mathbb{E}[\mathbf{f}(T)]\|_2\leq\frac{\sqrt{N}\alpha}{1-\|\mathbf{\Sigma}\|_2}\left(8\beta\max_{n=1,...,N}|\Omega_n|\max_{n=1,...,N}\|\mathbf{\widetilde{w}}_n\|_2+\gamma\sqrt{M}\right).
\end{align}
To summarize, we have proved the following theorem.
\begin{thm}
For Algorithm \ref{a2}, under Assumptions 1-5, we have:
\begin{align}
\limsup_{T\rightarrow\infty}\|\mathbb{E}[\mathbf{f}(T)]\|_2\leq\frac{\sqrt{N}\alpha}{1-\|\mathbf{\Sigma}\|_2}\left(8\beta\max_{n=1,...,N}|\Omega_n|\max_{n=1,...,N}\|\mathbf{\widetilde{w}}_n\|_2+\gamma\sqrt{M}\right),
\end{align}
where $\|\mathbf{\Sigma}\|_2<1$ is a constant related to the network structure and the error vector $\mathbf{f}(T)$ is defined as:
\begin{align}
\mathbf{f}(T)=
\left[
\begin{array}{c}
\mathbf{w}_1(T)-\mathbf{\widetilde{w}}_1\\
\vdots\\
\mathbf{w}_N(T)-\mathbf{\widetilde{w}}_N
\end{array}
\right]\label{convergence_subgradient}
\end{align}
\end{thm}

Two remarks regarding the result and assumptions of the theorem are in order, respectively.

\begin{rem}
From the R.H.S. of \eqref{convergence_subgradient}, it seems that the smaller the values of $\beta,\gamma$, the better the performance of Algorithm \ref{a2}. However, this is not true. The reason is that the R.H.S. of \eqref{convergence_subgradient} is only an upper bound on the error and when deriving this upper bound, the regularization terms are treated as biases and we want to bound their influences (similar remarks can be made to the analysis of most regularized adaptive algorithms, e.g., \cite{liu2014distributed,babadi2010sparls}). However, in practice, as we have argued, the facts that (i) neighbors have similar weights; and (2) weights are sparse; are important prior knowledge, which can enhance the performance of the algorithm. To embody this prior knowledge, we need to assign appropriate positive values to $\beta$ and $\gamma$. Another way to decrease the R.H.S. of \eqref{convergence_subgradient} is to decrease $\alpha$. However, this comes with a cost: smaller $\alpha$ leads to slower convergence speed. So, in practice, a positive value for $\alpha$ cannot be too small.
\end{rem}
\begin{rem}
Assumptions 1-3 are standard hypotheses to analyze the convergence behaviors of adaptive algorithms. Algorithm 5, i.e., small enough step size $\alpha$, is also common in gradient-descent-type algorithms such as LMS and subgradient method. The only nonconventional assumption is Assumption 4. In particular, Assumption 4 holds when $\lambda$ is sufficiently close to 1, which is the case (e.g., $\lambda=0.995,0.999$) for most applications of RLS \cite{Haykin:1996:AFT:230061,babadi2010sparls,liu2014distributed}. This is not surprising because in most applications the weight vectors vary fairly slowly (if the weight vectors change too drastically across time, then virtually no adaptive algorithms can track them well as the past data become useless and mere one shot data at current time are not enough for estimating the high dimensional weight vectors) and thus $\lambda$ is chosen very close to 1. Another interesting observation is that for networks with uniform degree, i.e., $\max_{n}|\Omega_n|=\min_{n}|\Omega_n|$, Assumption 4 is always satisfied.
\end{rem}

\section{Numerical Evaluation}

In this section, numerical simulations are conducted to verify the effectiveness of the proposed decentralized online ADMM algorithm (or ADMM algorithm in short), Algorithm 1, and the proposed decentralized online subgradient method (or subgradient method in short), Algorithm 2. The performance of the global offline optimization of problem \eqref{RLS} (global optimizor henceforth) is shown as a benchmark. For the sake of comparison, the performance of the distributed single task sparse RLS algorithm in \cite{liu2014distributed} (DSPARLS henceforth) is also presented to highlight the impact of multitask.

\begin{figure}
  \centering
  \includegraphics[scale=.2]{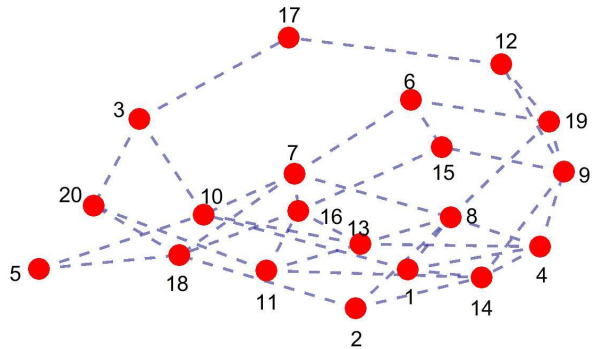}\\
  \caption{The network topology.}\label{network}
\end{figure}

We consider a network with $N=20$ nodes and 40 random edges so that the average node degree is 4. The network topology is illustrated in Fig. \ref{network}. The dimension of the input data is $M=20$. Each entry of the input data sequence $\{\mathbf{u}_n(t)\}$ is generated according to the uniform distribution over the interval $[0,1]$ independently. The noise sequence $\{e_n(t)\}$ is generated according to the uniform distribution on $[0,N_0]$ independently, where $N_0$ is a constant controlling the noisy level of the observations. To achieve sparsity, we let 18 entries (whose positions are randomly selected) of the true weight vectors $\widetilde{\mathbf{w}}_n(t)$ be zero. The two remaining entries $\widetilde{\mathbf{w}}_n^\text{part}(0)\in\mathbb{R}^2$ of the initial weight vectors are generated in a way that enforces similarity between neighbors. Specifically, we first generate $N$ i.i.d. two dimensional random vectors $\{\boldsymbol{\phi}_n\}_{n=1,...,N}$ uniformly distributed on $[0,1]^2$. Then, we solve the following optimization problem to obtain $\widetilde{\mathbf{w}}_n^\text{part}(0)$:
\begin{align}
\{\widetilde{\mathbf{w}}_n^\text{part}(0)\}_{n=1,...,N}=\arg\min_{\mathbf{w}_n\in\mathbb{R}^2,n=1,...,N}\sum_{n=1}^N\|\mathbf{w}_n-\boldsymbol{\phi}_n\|_2^2+\frac{1}{2}\sum_{n=1}^N\sum_{m\in\Omega_n}\|\mathbf{w}_n-\mathbf{w}_m\|_2^2,
\end{align}
which promotes similarity between neighbors and can be easily solved as the objective function is a convex quadratic function. To capture the slowly time-variant trait of the weight vectors, the increment from $\widetilde{\mathbf{w}}_n(t)$ to $\widetilde{\mathbf{w}}_n(t+1)$, i.e., $\widetilde{\mathbf{w}}_n(t+1)-\widetilde{\mathbf{w}}_n(t)$ is generated by uniform distribution on $[-0.5N_1,0.5N_1]$ independently across time and nodes, where $N_1$ is a constant controlling the varying rate of the weight vectors.

\begin{figure}
\renewcommand\figurename{\small Fig.}
\centering \vspace*{8pt} \setlength{\baselineskip}{10pt}

\subfigure[Learning curve of Scenario 1]{
\includegraphics[scale = 0.15]{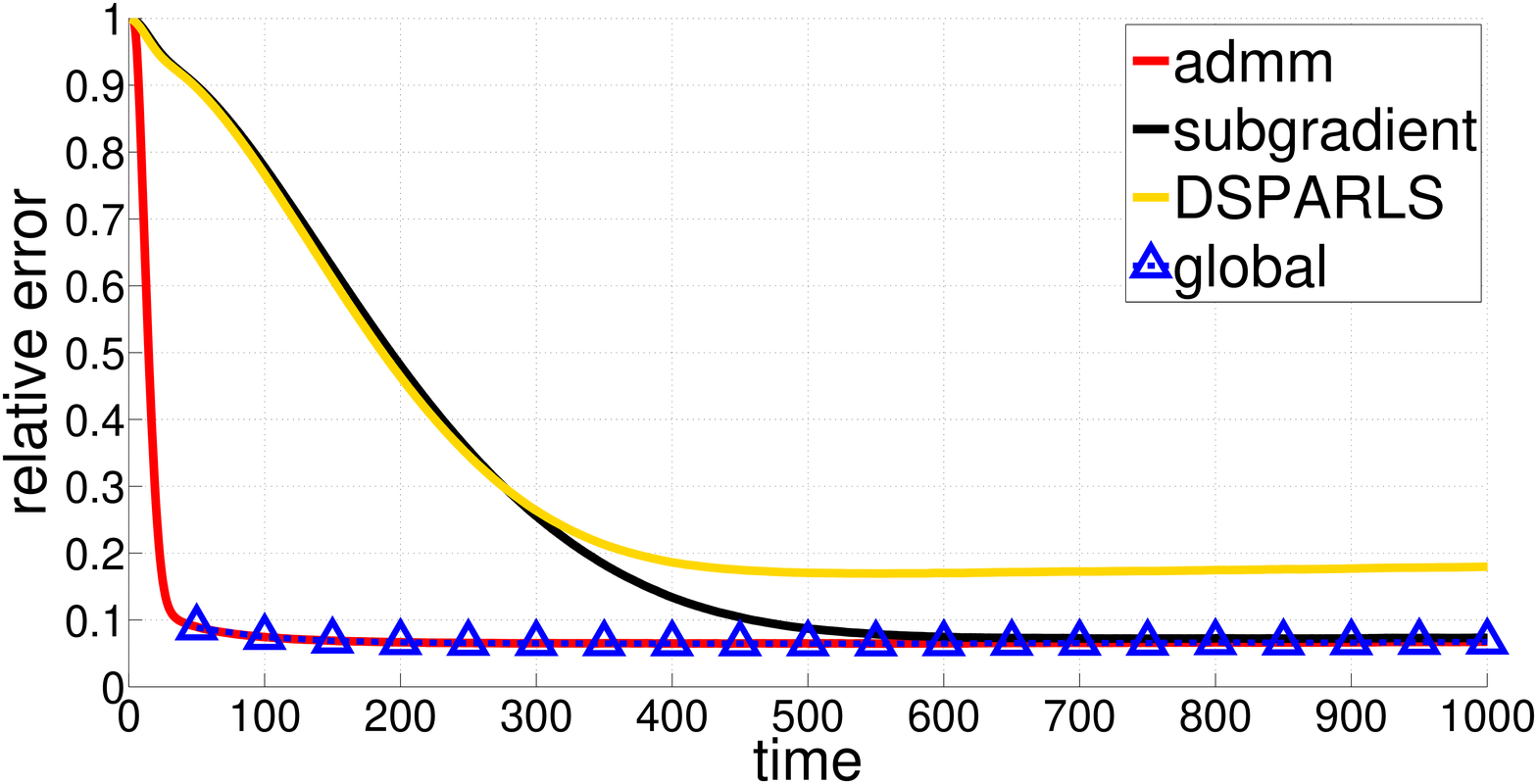}}
\subfigure[Learning curve of Scenario 2]{
\includegraphics[scale = 0.15]{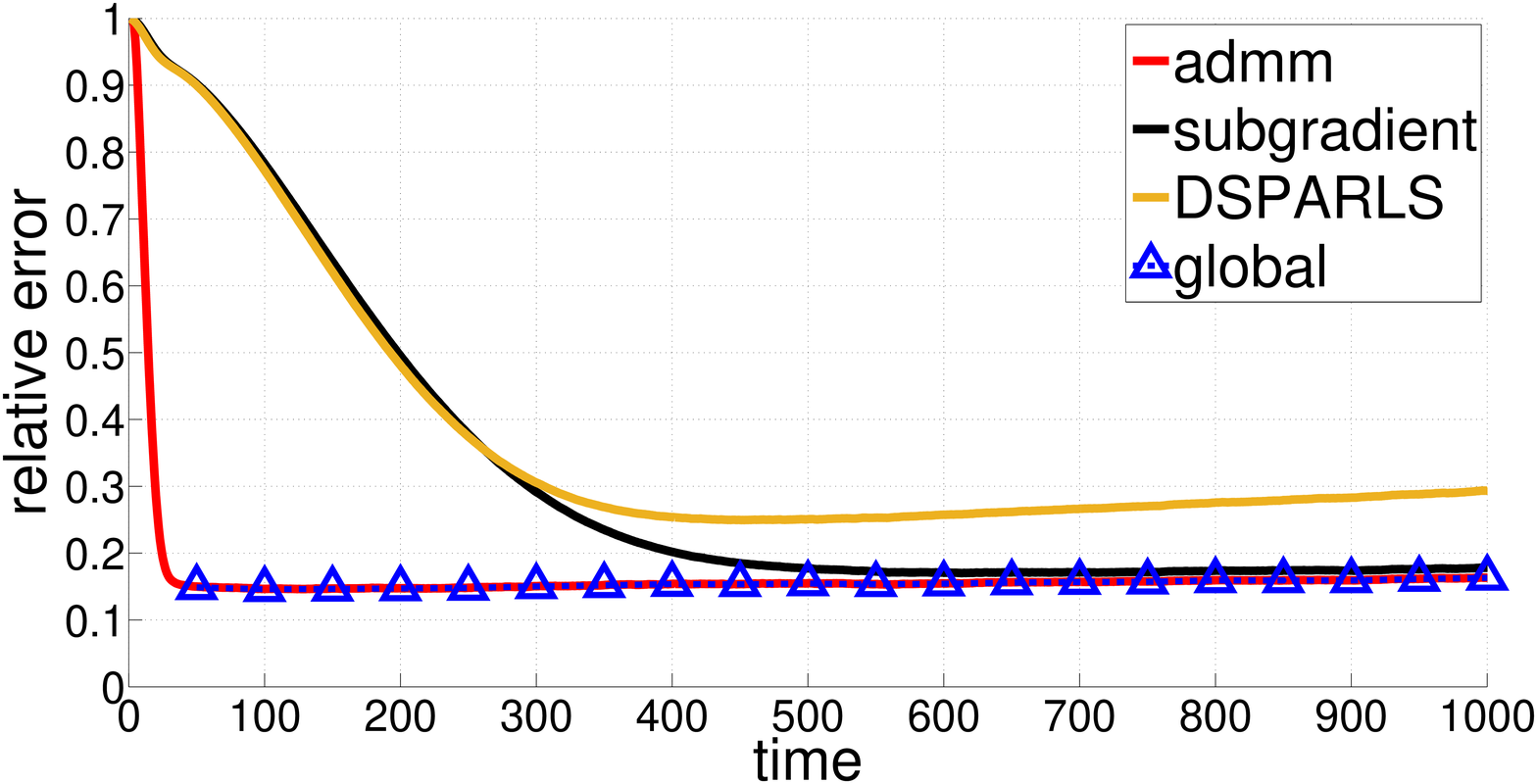}}
\caption{Learning curves of Scenario 1 and Scenario 2.}
\label{learning_curve}
\end{figure}

Now, we choose the regularization parameters and forgetting factor as $\beta=\gamma=1,\gamma=0.995$. We consider two scenarios in terms of the noise level $N_1$ and the varying rate of weight vectors $N_2$. In Scenario 1, $N_1=0.1,N_2=0.02$ while in Scenario 2, $N_1=0.3,N_2=0.05$. The latter scenario has noisier observations and weight vectors which vary faster. Thus, the weight vectors of Scenario 2 are more difficult to track than those of Scenario 1. For the proposed ADMM algorithm, the proposed subgradient method, the DSPARLS algorithm in \cite{liu2014distributed} and the global optimizor, we plot the relative errors (defined to be $\|\mathbf{w}(t)-\widetilde{\mathbf{w}}(t)\|_2/\|\widetilde{\mathbf{w}}(t)\|_2$, where $\mathbf{w}(t)$ and $\widetilde{\mathbf{w}}(t)$ are concatenations of $\mathbf{w}_n(t)$ and $\widetilde{\mathbf{w}}_n(t)$ of all nodes, respectively) as functions of time indices, i.e., the learning curves, under Scenario 1 (Fig. \ref{learning_curve}-(a)) and Scenario 2 (Fig. \ref{learning_curve}-(b)), respectively. Each learning curve is the average of 300 independent trials. Several interesting observations can be made from Fig. \ref{learning_curve}. First, the relative errors of both the proposed ADMM algorithm and the proposed subgradient method can converge to that of the global optimizor, i.e., the performance benchmark, as the observation data accumulate. On the contrary, the relative error of DSPARLS does not converge to that of the global optimizor. This highlights the effectiveness of the two proposed algorithms when tracking multitask weight vectors, which cannot be tracked well by existing method (DSPARLS in this case) for the single task situation. Second, comparisons between the learning curves of the proposed two algorithms indicate that the proposed ADMM algorithm needs much fewer observations, or equivalently much less time (about 100 time units), to track the weight vectors accurately than the proposed subgradient method does (about 600 time units). This is not surprising as dual domain methods generally converge faster than primal domain methods in the literature of optimization theory \cite{ling2015dlm,mokhtari2016dqm}. However, the advantage of the proposed ADMM algorithm in convergence speed comes at the cost of higher computational overhead per time unit than the proposed subgradient method. This accuracy-complexity tradeoff makes the proposed two algorithms appropriate for different applications depending on the computational capability of devices and needed tracking accuracy. Third, as one expects, Scenario 1 has better tracking performance than Scenario 2: the ultimate relative error of the proposed algorithms in Scenario 1 is about 0.067 while that of Scenario 2 is about 0.17. So, higher noise level and faster varying speed of the weight vectors do result in lower tracking accuracy.

\begin{figure}
\renewcommand\figurename{\small Fig.}
\centering \vspace*{8pt} \setlength{\baselineskip}{10pt}

\subfigure[Relative errors of each node at time 200 in Scenario 1]{
\includegraphics[scale = 0.15]{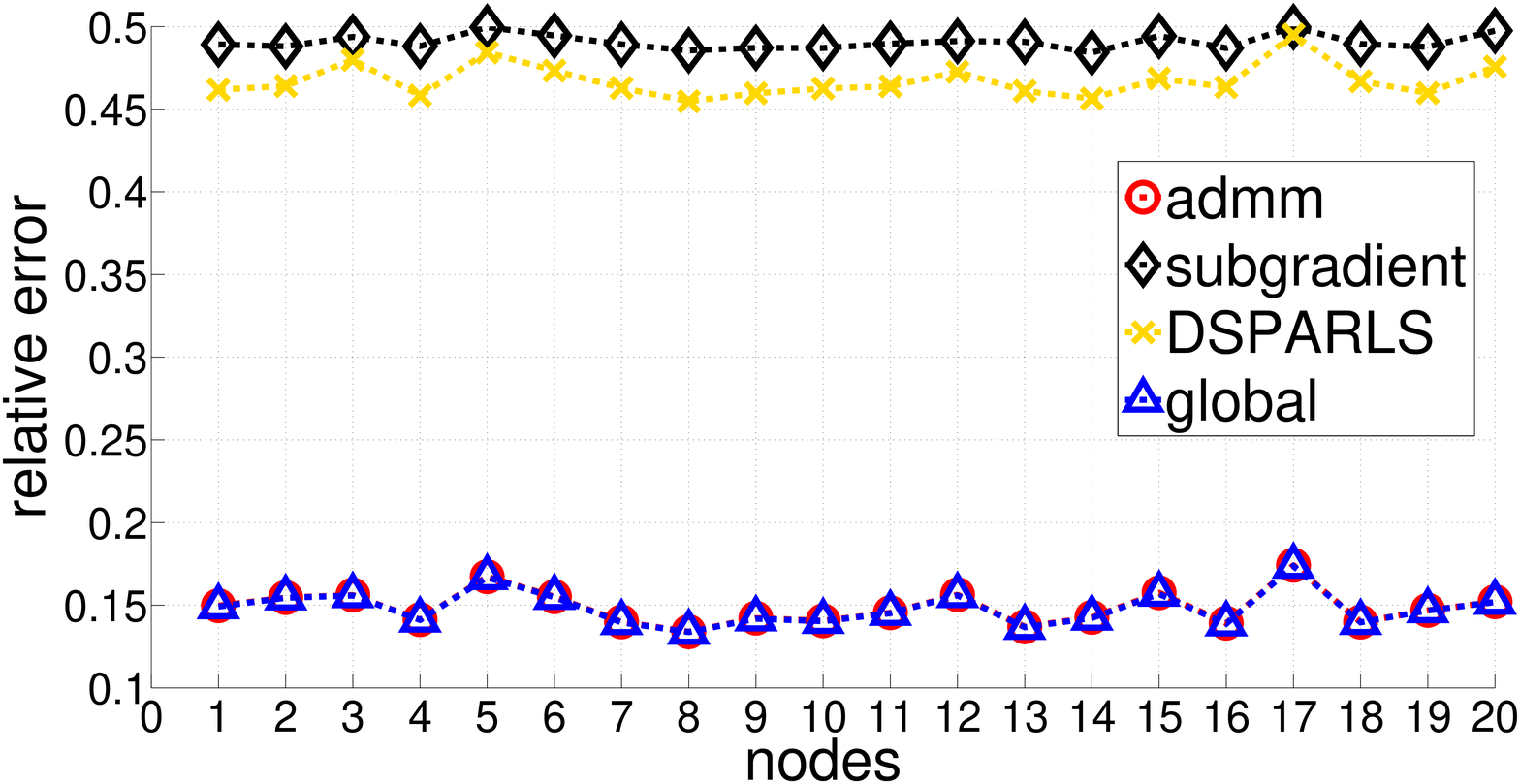}}
\subfigure[Relative errors of each node at time 500 in Scenario 1]{
\includegraphics[scale = 0.15]{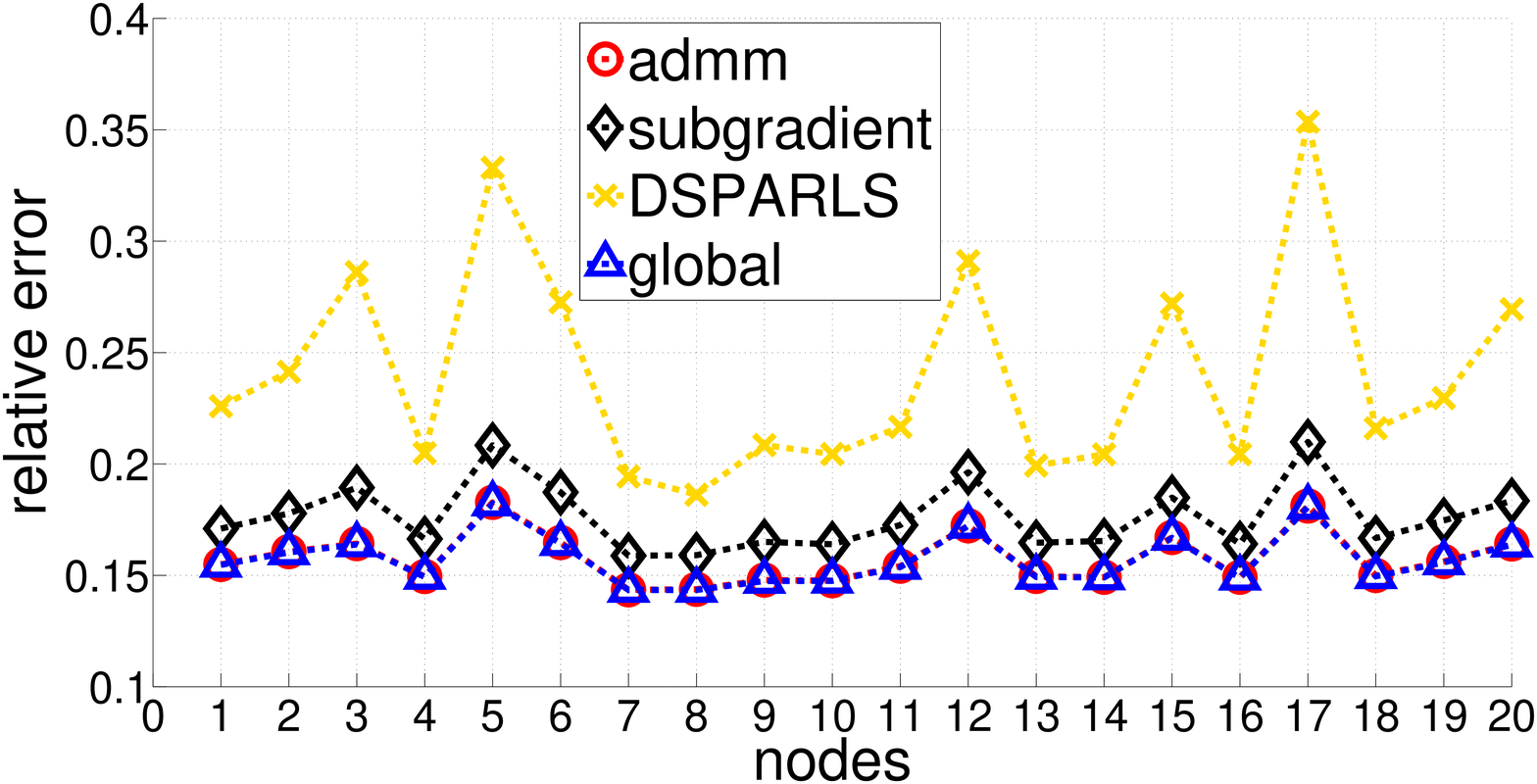}}
\subfigure[Relative errors of each node at time 200 in Scenario 2]{
\includegraphics[scale = 0.15]{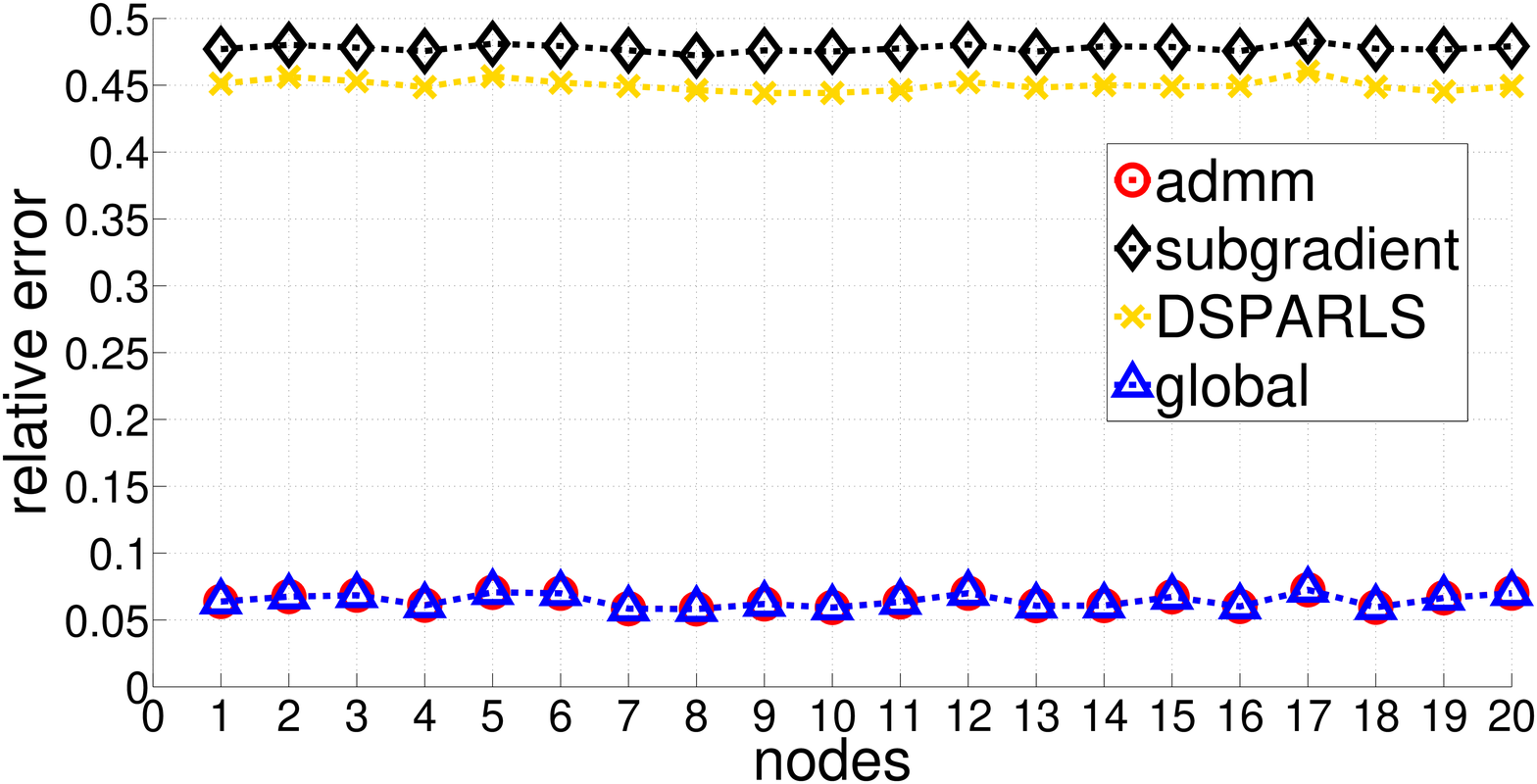}}
\subfigure[Relative errors of each node at time 500 in Scenario 2]{
\includegraphics[scale = 0.15]{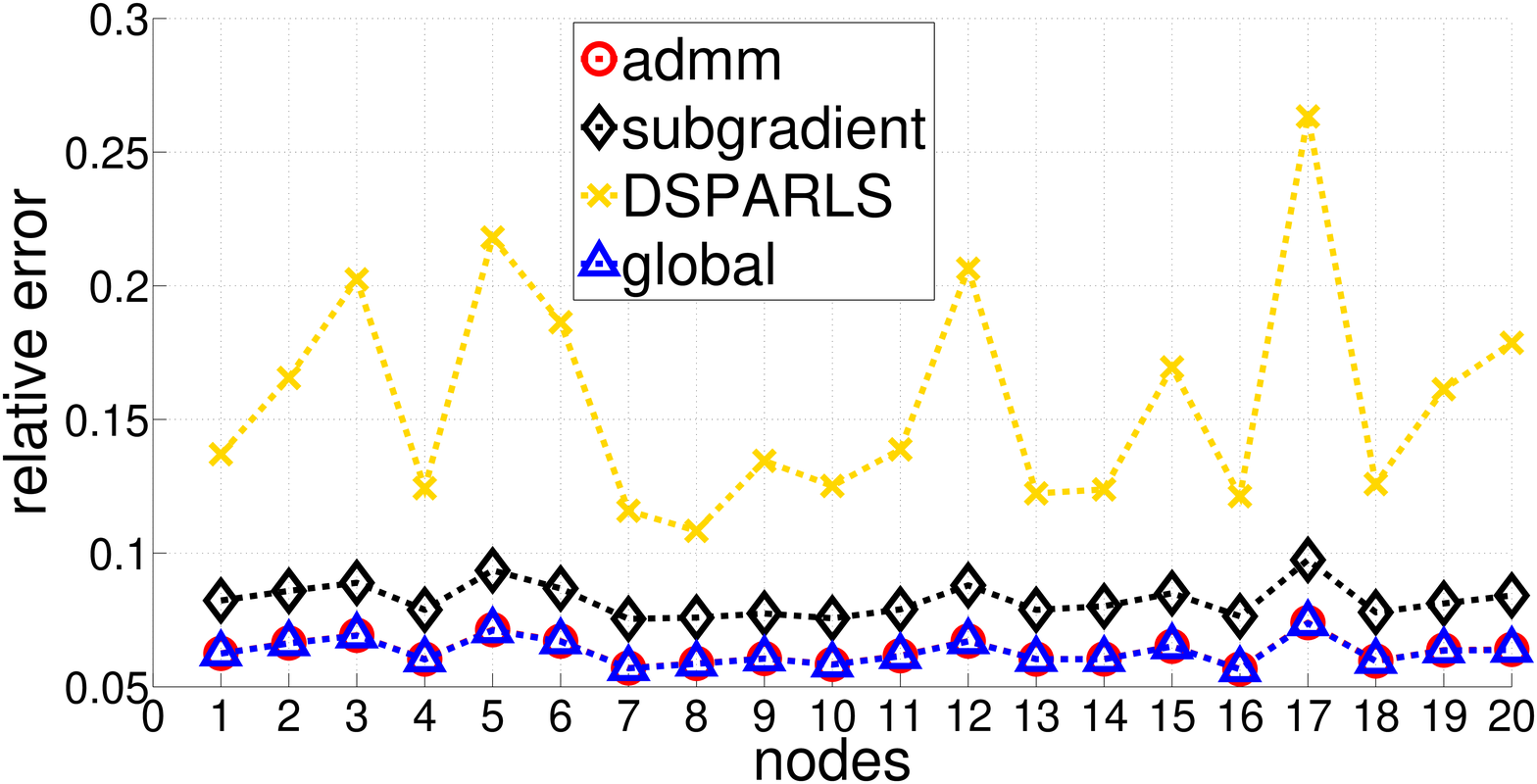}}

\caption{Relative tracking errors of each node.}
\label{nodes}
\end{figure}

Next, we investigate the tracking performance of each individual node. To this end, in Fig. \ref{nodes}, we show the relative tracking errors of each node at time $200$ and $500$ under Scenario 1 and Scenario 2. Several remarks are in order. First, we note that in all four cases of Fig. \ref{nodes}, the red curve (the proposed ADMM algorithm) and the blue curve (the global optimizor) coincide precisely for every node. This further confirms the previous observation from Fig. \ref{learning_curve} that the performance of the proposed ADMM can converge to that of the global optimizor quickly. Second, the proposed subgradient method, though performs poorly at time 200, has relative errors close to those of the global optimizor at time 500. This suggests that the proposed subgradient method eventually has performance close to the benchmark (the global optimizor). But this good performance necessitates longer time (or equivalently more data) compared to the proposed ADMM algorithm. Third, the performance of the single task learning algorithm DSPARLS never converge to that of the global optimizor. In particular, from Fig. \ref{nodes}-(b)(d), the performance of DSPARLS is worst at node 5 and node 17. Recall the network topology in Fig. \ref{network} and we see that these two nodes are loosely connected to other nodes, i.e., their degree is low. Thus, the weight vectors at these two nodes can potentially deviate far from the weight vectors at the rest of the nodes and thus violate the single task assumption of DSPARLS the most among all nodes. This partially explains the poor performance of DSPARLS at nodes 5 and 17.

\begin{figure}
\renewcommand\figurename{\small Fig.}
\centering \vspace*{8pt} \setlength{\baselineskip}{10pt}

\subfigure[Number of successful trials for different $\beta$]{
\includegraphics[scale = 0.15]{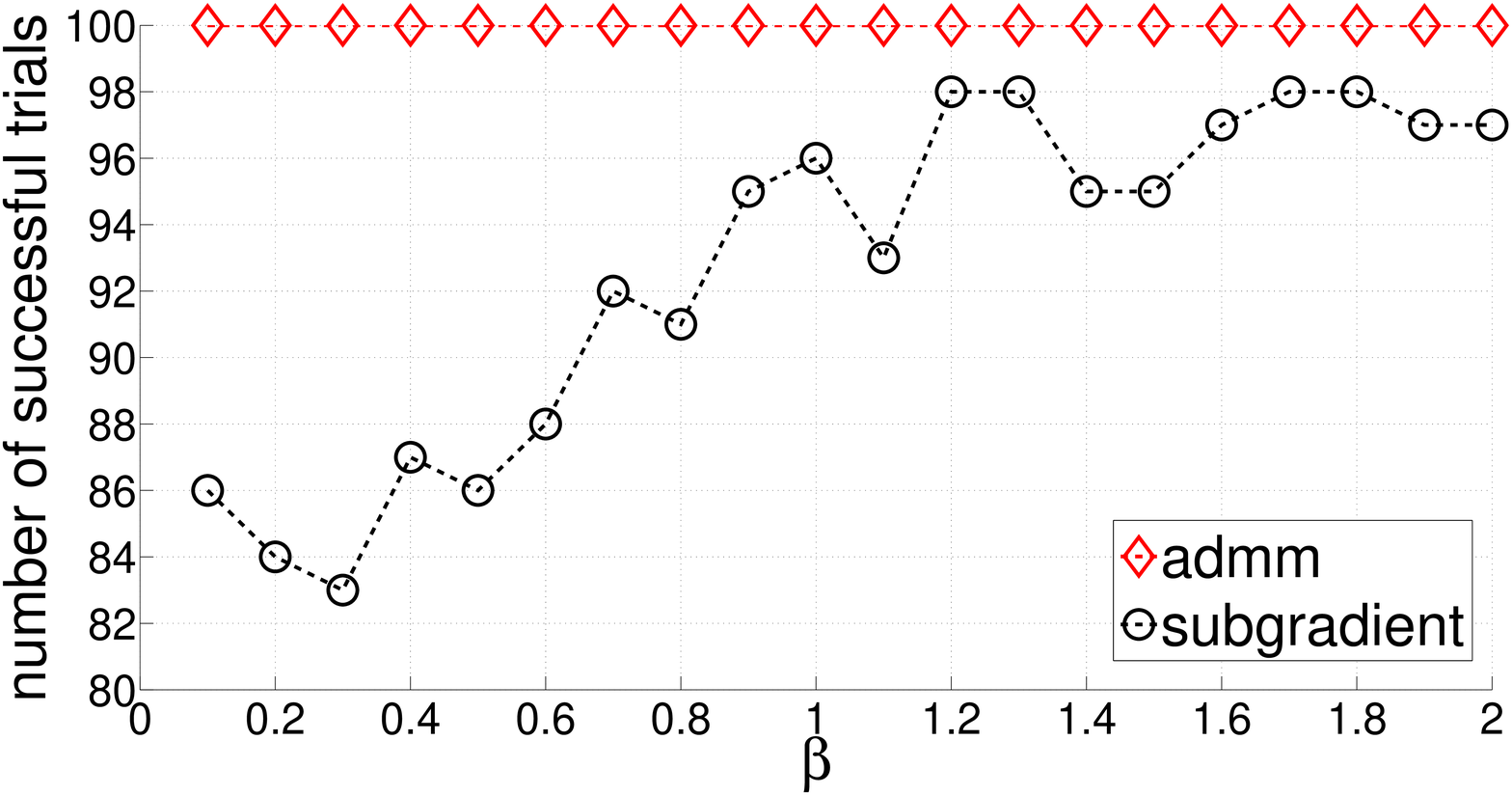}}
\subfigure[Average time needed to reach success among successful trials for different $\beta$]{
\includegraphics[scale = 0.15]{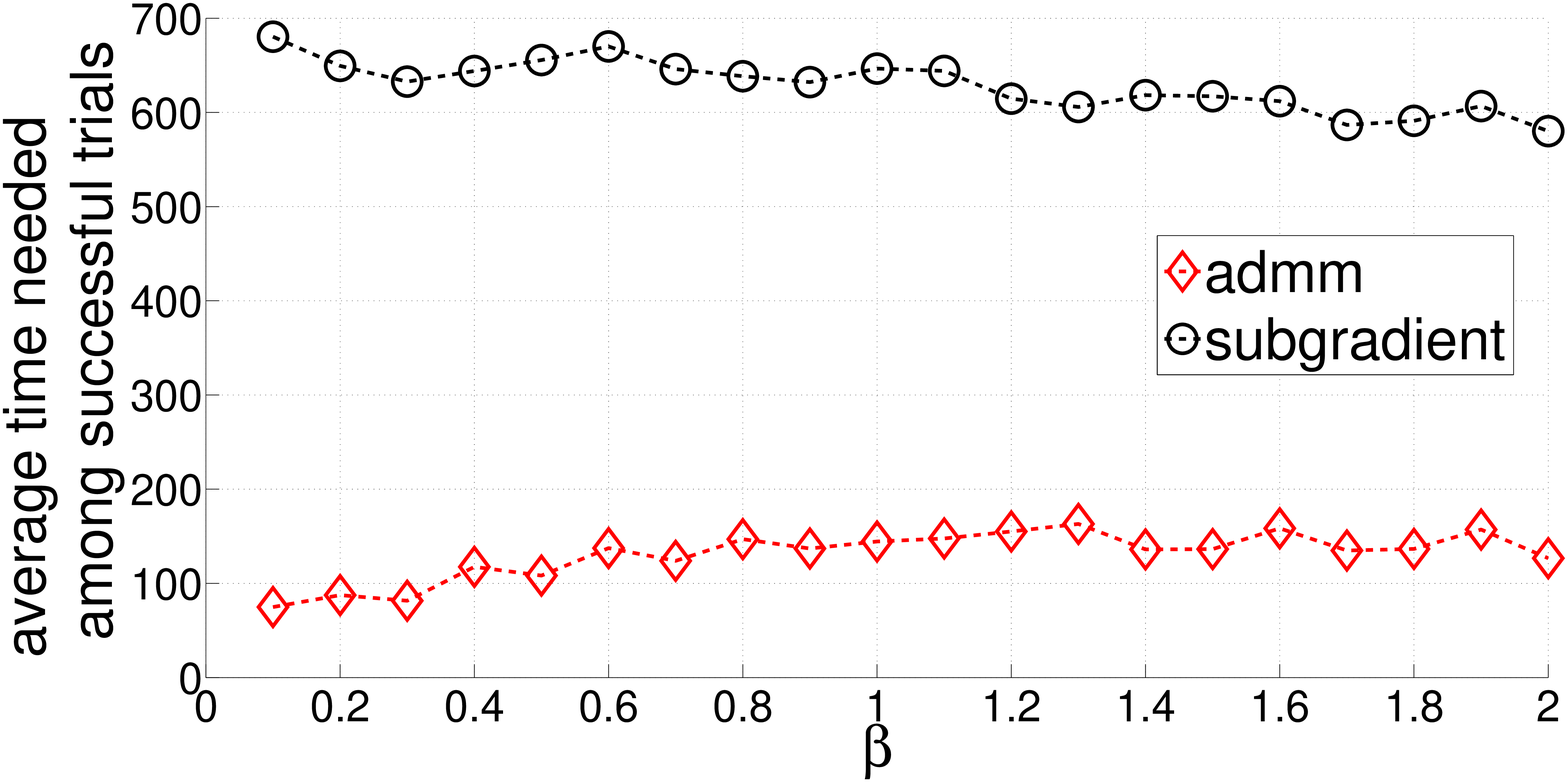}}

\subfigure[Number of successful trials for different $\gamma$]{
\includegraphics[scale = 0.15]{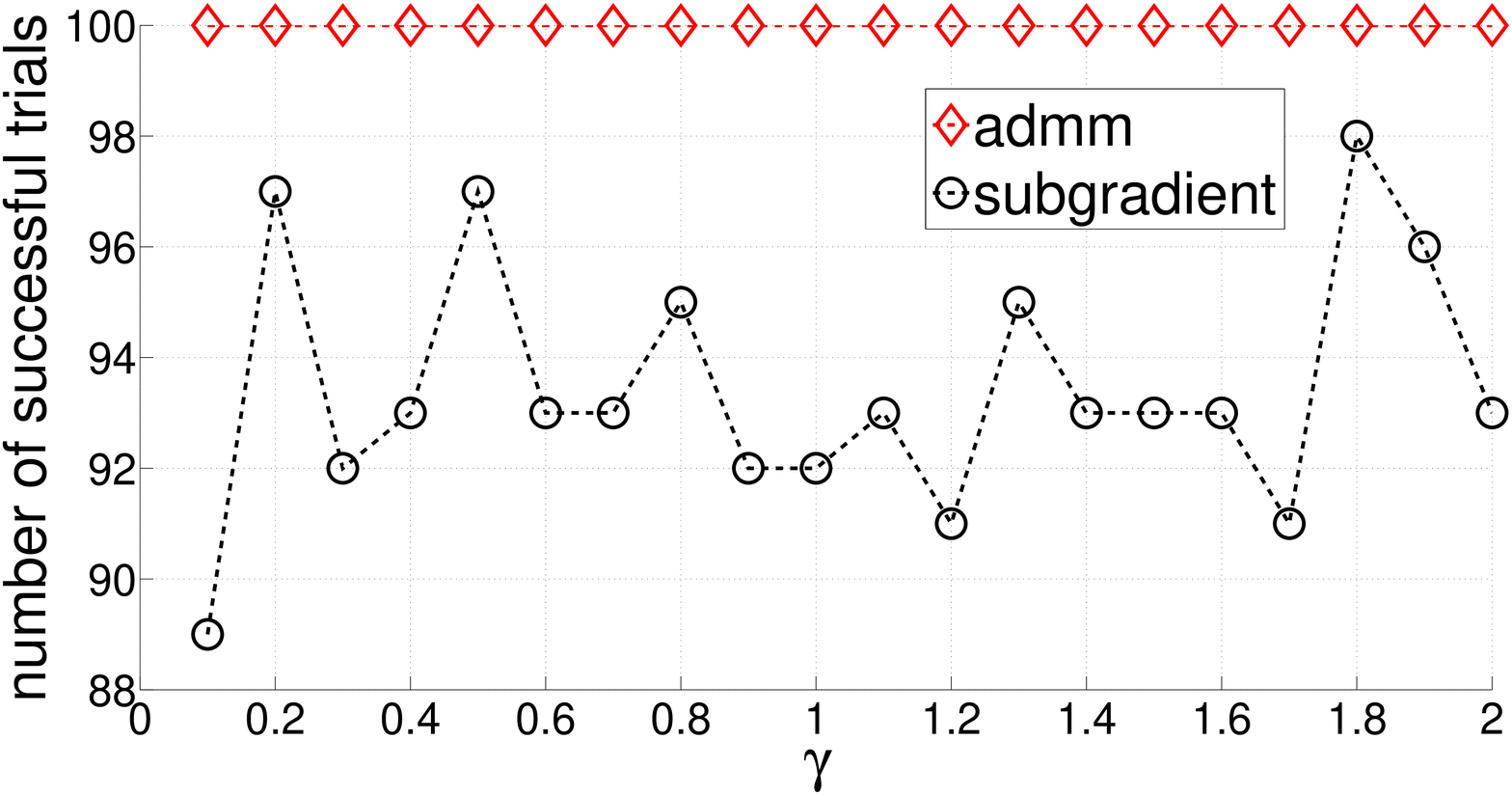}}
\subfigure[Average time needed to reach success among successful trials for different $\gamma$]{
\includegraphics[scale = 0.15]{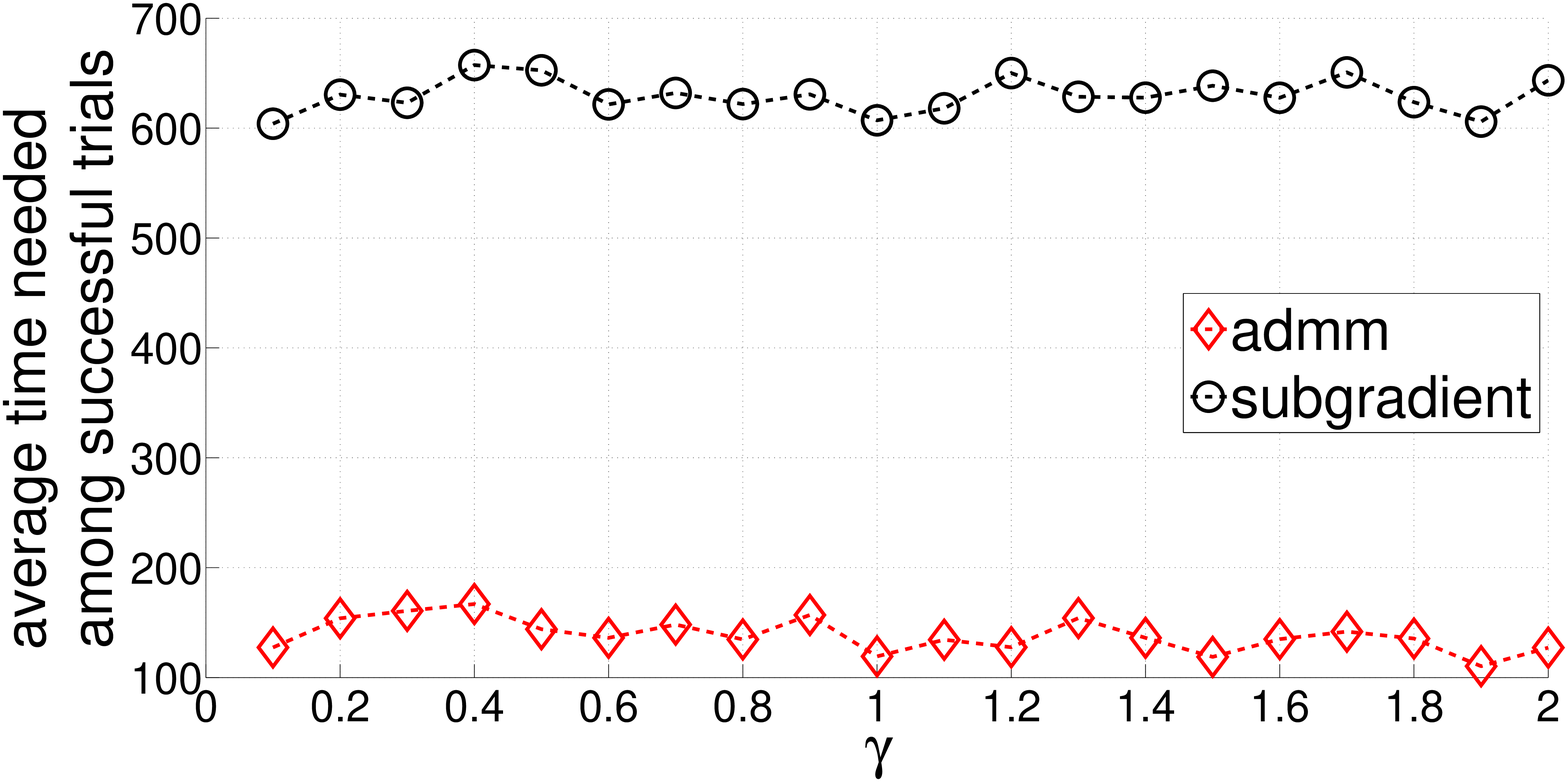}}

\subfigure[Number of successful trials for different $\lambda$]{
\includegraphics[scale = 0.15]{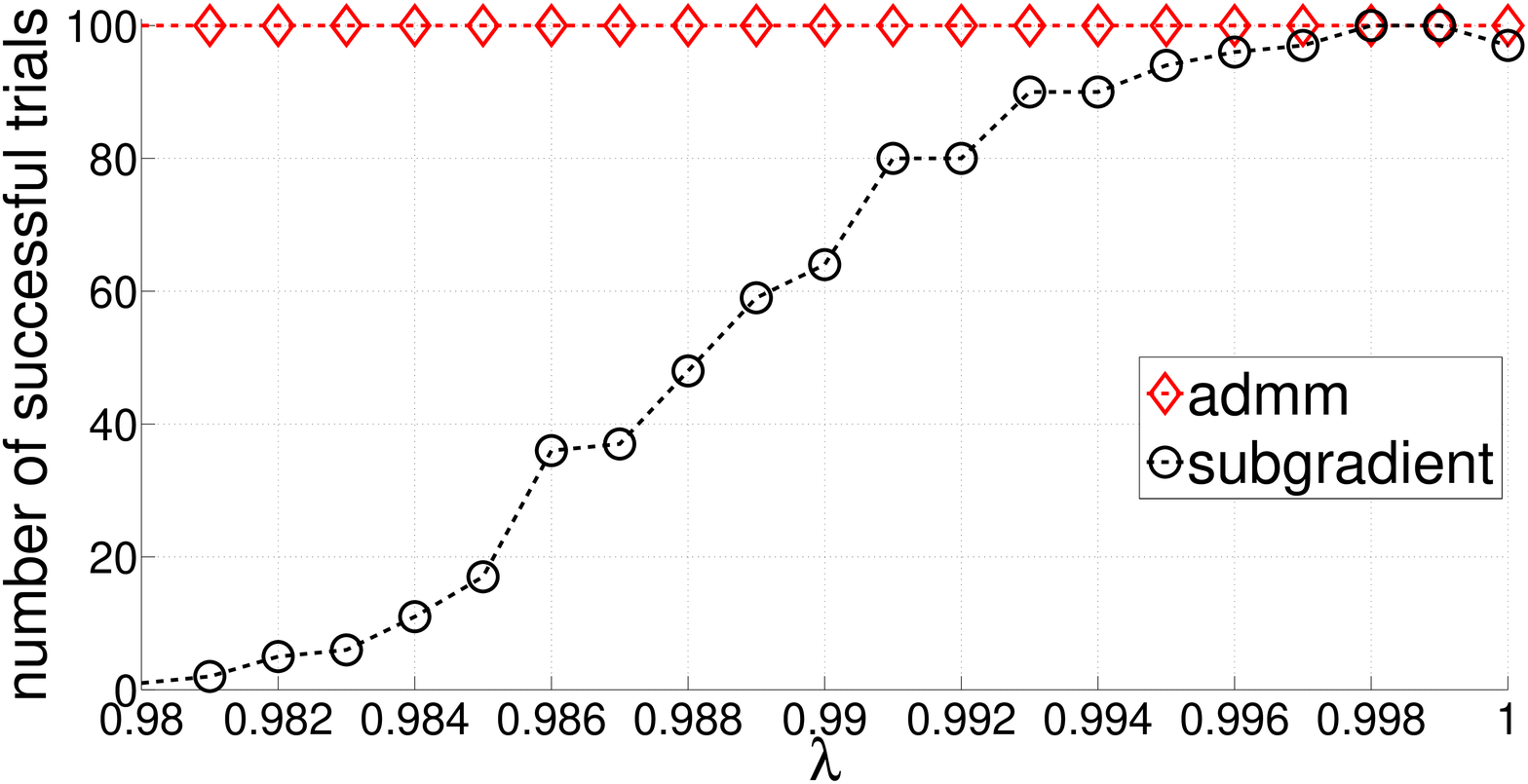}}
\subfigure[Average time needed to reach success among successful trials for different $\lambda$]{
\includegraphics[scale = 0.15]{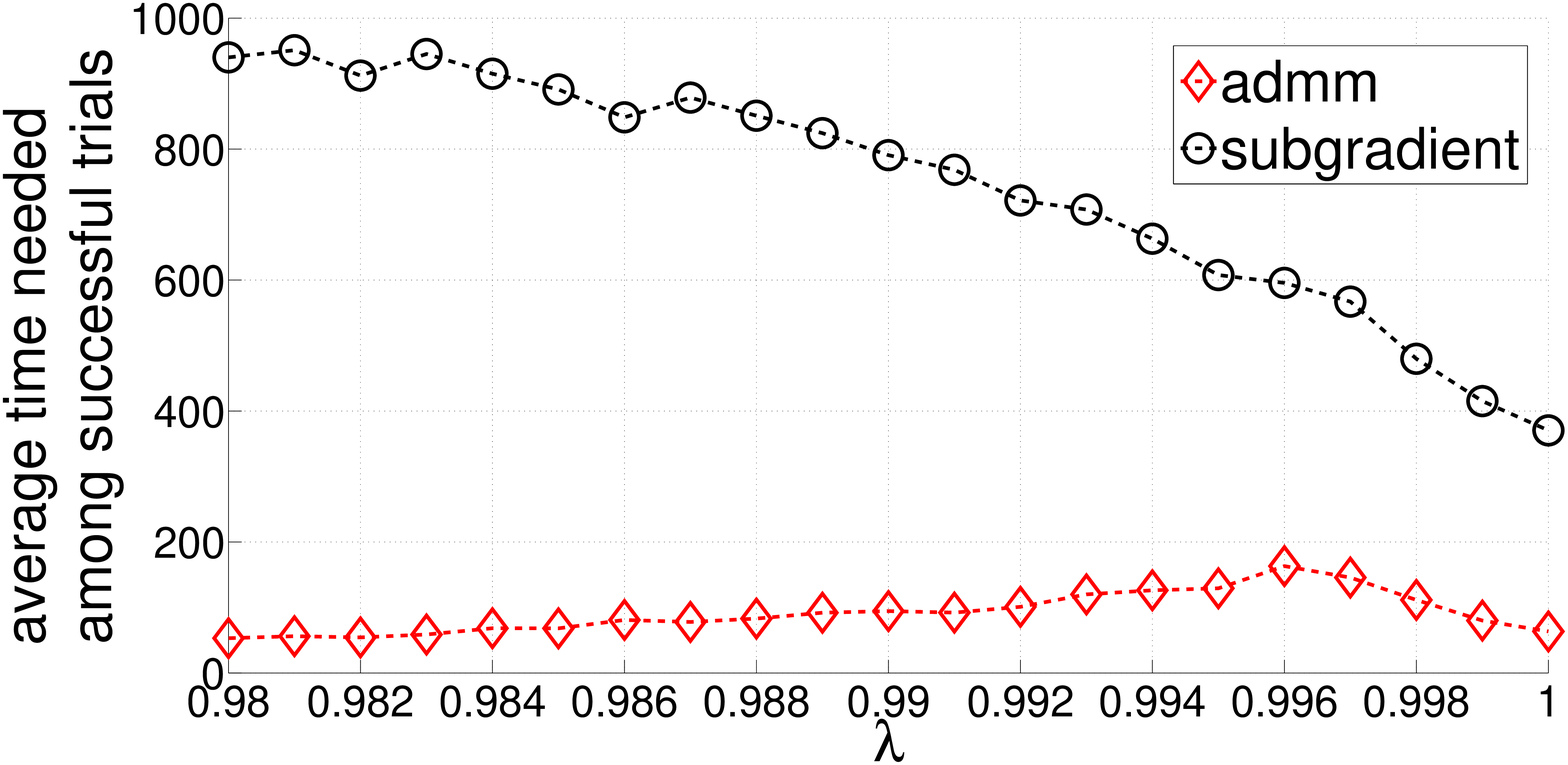}}

\caption{Number of successful trials and the average time needed to reach success among successful trials.}
\label{parameter}
\end{figure}

Previous experiments indicate that the proposed ADMM algorithm possesses faster and more accurate tracking performance than the proposed subgradient method. Next, we conduct a more thorough performance comparison between the proposed two algorithms for different regularization parameters $\beta,\gamma$ and different forgetting factors $\lambda$. We note that the global optimizor usually converges well before time $1000$ and we denote the steady relative error of the global optimizor at time $1000$ as $\check{e}$. We say a simulation trial of an algorithm (either the proposed ADMM algorithm or the proposed subgradient method) is \emph{successful} if, before time 1000, there exists a time window (i.e., interval) of length 20 over which the average relative error of the algorithm is lower than $1.1\check{e}$. The basic parameter setup is $\gamma=\beta=1,\lambda=0.995$. In each of the subfigures in Fig. \ref{parameter}-(a)(c)(e), we vary one parameter while keep the remaining two parameters the same as the basic setup. For each parameter setup, we conduct 100 independent trials and plot the number of successful trials in Fig. \ref{parameter}-(a)(c)(e). We observe that (i) the proposed ADMM algorithm is always successful, i.e., it can always converge to the steady performance of the global optimizor; (ii) the proposed subgradient method is successful in most trials as long as the forgetting factor $\lambda$ is sufficiently close to 1, which is the case in most applications (e.g., $\lambda=0.995$) as the weight vectors are varying very slowly and a large $\lambda$ is needed for tracking them. Note that the large $\lambda$ is reminiscent of the Assumption 4 in the performance analysis of the proposed subgradient method. See Remark 3 for a justification of large $\lambda$. Moreover, we investigate the average time needed to reach success (defined to be the middle point of the first time window over which the average relative error is less than $1.1\check{e}$) among successful trials. The results are shown in Fig. \ref{parameter}-(b)(d)(f). We remark that the proposed ADMM mostly needs no more than 150 time units to be successful, i.e., be close to the steady performance of the global optimizor, while it takes the proposed subgradient method a much longer time (around 600 time units) to be successful. This further confirms our previous assertion that the proposed ADMM algorithm possesses faster tracking performance than the proposed subgradient method.

\section{Conclusion}

In this paper, we study the decentralized sparse multitask RLS problem. We first propose a decentralized online ADMM algorithm for the formulated RLS problem. We simplify the algorithm so that each iteration of the ADMM algorithm consists of simple closed form computations and each node only needs to store and update one $M\times M$ matrix and six $M$ dimensional vectors. We show that the gap between the iterations of the proposed ADMM algorithm and the optimal points of the formulated RLS problem will converge to zero. Moreover, to further reduce the computational complexity, we propose a decentralized online subgradient method. We show that the iterations of the proposed subgradient method will track the true weight vectors with errors upper bounded by some constant related to the network topology and algorithm parameters. Compared with the ADMM algorithm, the subgradient method enjoys lower computational complexity at the expense of slower convergence speed. The proposed algorithms are corroborated by numerical experiments, which highlight the effectiveness of the proposed algorithms and the accuracy-complexity tradeoff between the proposed two algorithms.

\bibliography{mybib}{}
\bibliographystyle{ieeetr}

\end{document}